\documentclass[conference]{IEEEtran}
\IEEEoverridecommandlockouts
\usepackage{amsmath}
\usepackage{graphicx, subfigure}
\usepackage{amssymb}
\usepackage{cite}
\usepackage{url}
\usepackage[ruled,vlined,lined,ruled,linesnumbered]{algorithm2e}
\usepackage{relsize}
\usepackage[nodisplayskipstretch]{setspace}
\usepackage{longtable}
\usepackage{booktabs}
\usepackage{multirow}
\usepackage[subnum]{cases}

\newtheorem{theorem}{Theorem}
\newtheorem{lemma}{Lemma}
\newtheorem{remark}{Remark}

\newtheorem{proposition}{Proposition}

\begin{document}

\title{ A Model Predictive Control Approach for Low-Complexity Electric Vehicle Charging Scheduling: Optimality and Scalability }

\author{\IEEEauthorblockN{Wanrong Tang ~\IEEEmembership{Student Member,~IEEE} and Ying Jun (Angela) Zhang,~\IEEEmembership{Senior Member,~IEEE}}

\IEEEauthorblockA{Department of Information Engineering, The Chinese University of Hong Kong\\
Shatin, New Territories, Hong Kong
}

\thanks{This work was presented in part as \cite{tang2014onlineMSP}. This work was supported in part by the National Basic Research Program (973 program Program number 2013CB336701), and three grants from the Research Grants Council of Hong Kong under General Research Funding (Project number 2150828 and 2150876) and Theme-Based Research Scheme (Project number T23-407/13-N). }

\thanks{W.~Tang is with the Department of Information Engineering, The Chinese University of Hong Kong, Shatin, New Territories, Hong Kong (Email: twr011@ie.cuhk.edu.hk).}

\thanks{Y. J. Zhang is with the Department of Information Engineering, The Chinese University of Hong Kong. She is also with Shenzhen Research Institute, The Chinese University of Hong Kong (Email: yjzhang@ie.cuhk.edu.hk).}

\vspace{-1cm}
}

\maketitle

\begin{abstract}

With the increasing adoption of plug-in electric vehicles (PEVs), it is critical to develop efficient charging coordination mechanisms that minimize the cost and impact of PEV integration to the power grid.
In this paper, we consider the optimal PEV charging scheduling,
where the non-causal information about future PEV arrivals is not known in advance, but its statistical information can be estimated.
This leads to an ``online'' charging scheduling problem that is naturally formulated as a finite-horizon dynamic programming with continuous state space and action space.
To avoid the prohibitively high complexity of solving such a dynamic programming problem, we provide a Model Predictive Control (MPC) based algorithm with computational complexity $O(T^3)$, where $T$ is the total number of time stages.
We rigorously analyze the performance gap between the near-optimal solution of the MPC-based approach and the optimal solution for any distributions of exogenous random variables.
Furthermore, our rigorous analysis shows that when the random process describing the arrival of charging demands is first-order periodic, the complexity of proposed algorithm can be reduced to $O(1)$, which is independent of $T$.
Extensive simulations show that the proposed online algorithm performs very closely to the optimal online algorithm.
The performance gap is smaller than $0.4\%$ in most cases.

\end{abstract}

\section{Introduction}
\subsection{Background and Contributions}\label{subsec:background and contribution}
The massive deployment of PEVs imposes great challenges to smart power grid, such as voltage deviation, increased power losses, and higher peak load demands. It is critical to design PEV charging mechanisms that minimize the cost and impact of PEV integration.
Previously, PEV charging coordination has been extensively studied to  minimize power loss, minimize load variance, or minimize charging cost, etc \cite{sortomme2011coordinated,tang2014online,ma2013decentralized,he2012optimal,gan2012optimal}.
Ideally, the load demand can be flattened as much as possible if the information about future charging demand is known non-causally when calculating the charging schedule.
However, in practice, a PEV charging station only knows the load demand of the PEVs that have arrived, but not that of the PEVs coming in the future.
Fortunately, the statistical information of the future charging demands can often be acquired through historic data, which benefits the control of the PEV charging scheduling in practical scenarios.

In this paper, we consider the optimal PEV charging scheduling, assuming that the future charging demand is not known a priori, but its statistical information can be estimated.
In particular, we define the cost of PEV charging as a general strictly convex increasing function of the instantaneous load demand. Minimizing such a cost leads to a flattened load demand, which is highly desirable for many reasons \cite{sortomme2011coordinated,tang2014online,ma2013decentralized,he2012optimal,gan2012optimal}.
The online PEV charging scheduling problem is formulated as a finite-horizon dynamic programming problem with continuous state space and action space.
To avoid the prohibitively high complexity of solving such a dynamic programming problem,
we provide a Model Predictive Control (MPC) approach to obtain a near-optimal solution.
Instead of adopting the generic convex optimization algorithms to solve the problem, we propose an algorithm with computational complexity $O(T^3)$ by exploring the load flattening feature of the solution, where $T$ is the total number of time stages.
We rigorously analyze the performance gap between the near-optimal solution of the MPC-based approach and the optimal solution, and the result applies to any distributions of exogenous random variables.
Specially, the performance gap is evaluated by the Value of the Stochastic Solution (VSS), which represents the gap between the solution of the approximate approach and that of dynamic programming problem \cite{birge1997introduction,defourny2011multistage,maggioni2012analyzing}.
Furthermore, our analysis shows that when the random process describing the arrival of charging demands is first-order periodic, the complexity of proposed algorithm can be reduced to $O(1)$, which is independent of $T$.
Extensive simulations show that the proposed algorithm performs very closely to the optimal solution. The performance gap is smaller than $0.4\%$ in most cases.
As such, the proposed algorithm is very appealing for practical implementation due to its scalable computational complexity and close to optimal performance.

The rest of the paper is organized as follows.
A review of the related work on the PEV charging scheduling with uncertain load demand is presented in Section \ref{subsec:related work}.
We introduce the problem formulations of both offline and online PEV charging problem in Section \ref{sec:problem formulation}.
In Section \ref{sec:online algorithm}, we propose a MPC based online algorithm and analyze its performance gap. The $O(1)$-complexity algorithm is given when the arrival process is first-order periodic  in Section \ref{sec:case_study}. Simulation results are presented in Section \ref{sec:simulation}. Finally, the paper is concluded in Section \ref{sec:conclusions}.

\subsection{Related Work}\label{subsec:related work}
The works on the PEV charging scheduling with uncertain PEV load demand include both simulation-based evaluations \cite{leou2014stochastic,rao2014smart} and theoretical performance guarantees \cite{chen2014distributional,gan2013realtime,bansal2014plug,li2012modeling,zhang2014charging}.
Meanwhile, MPC is one of most commonly approaches for which has been widely adopted in recent studies \cite{rao2014smart,chen2014distributional,gan2013realtime,bansal2014plug}.
\cite{rao2014smart} leverages the MPC based method to design a dynamic charging and driving cost control scheme.
Both \cite{chen2014distributional} and \cite{gan2013realtime} apply MPC algorithms to minimize the load variation.
\cite{bansal2014plug} proposes a plug and play MPC approach to minimize the voltage fluctuations by assuming that the load demand is time-periodic.
Compared to \cite{rao2014smart,chen2014distributional,gan2013realtime,bansal2014plug}, in this paper we analyze the performance gap between the solution of MPC approach and the optimal solution regardless of the distribution of the load demand.
Besides, we provide a more scalable algorithm with $O(1)-$ complexity as well as the optimality analysis for the case when the load demand is first-order periodic.
Additionally, the objective functions in \cite{rao2014smart,chen2014distributional,gan2013realtime,bansal2014plug} are quadratic forms of load demand. Whereas in this paper, the objective function is a general strictly convex increasing function which reflects both the charging cost and the load variance.

As to the amount of information needed, the EV charging scheduling algorithms in \cite{li2012modeling} and \cite{zhang2014charging} require the probability distribution of the random PEV arrival process. In contrast, the proposed algorithm in this paper only requires the first-order moment, i.e., the expected values of the random demand patterns. In practice, it is a lot easier to obtain the expected values of random process than to obtain the probability distribution of a random process.
Convex cost functions are also considered in \cite{koutsopoulos2011Optimal} and \cite{huang2012optimal}. Both of them devise the online algorithms for battery charging control problems, where there is no charging deadline for the battery. The PEV charging scheduling problem in this paper differs from stationary battery charging in that each PEV has a demand to be satisfied before a certain deadline.

\section{Problem Formulation}\label{sec:problem formulation}
We consider the PEV charging scheduling problem, where PEVs arrive at the charging station at random instants with random charging demands that must be fulfilled before a random departure time.

\subsection{Optimal Offline PEV Charging Problem}\label{subsec:offline problem}
For the ease of understanding, we first introduce an ideal case, where all the non-causal information of based load and the PEVs, including the arrival times, departure times and charging demands are known to the charging station before the system time.
The entire system time is divided into $T$ equal-length time slots.
Let $\mathcal{N}$ denote the set of PEVs that arrive during the system time.
Notice that for a given time slot number $T$, $\mathcal{N}$ is itself a random set due to the random arrival of PEVs.
We denote by $\mathcal{I}(t)$ the set of PEVs that are in the charging station during slot $t$. Denote by $t_i^{(s)}$ and $t_i^{(e)}$ the arrival time slot and the departure time slot of PEV $i$, respectively. $d_i$ denotes the charging demand that PEV $i$ requires. The charging station needs to decide, the charging rate $x_{it}, \forall i \in \mathcal{I}(t)$. To satisfy the demand $d_i$, $x_{it}$ must satisfy
$\sum_{t = t_i^{(s)}}^{t_i^{(e)}} x_{it} = d_i.$
Let $s_t$ be the total charging rate of time slot $t$, i.e.,
\begin{equation}\small
s_t = \sum_{i \in \mathcal{I}(t)} x_{it}, \forall t = 1, 2, \cdots, T,
\end{equation}
which is also called charging load at time $t$.
The total load consists of both the charging load and the inelastic base load in the same location. The base load, denoted by $l_t$, represents the load of other electricity consumptions at time $t$ except for PEV charging. Then, the total load
at time $t$ is given by $\sum_{i \in \mathcal{I}_t}x_{it}+l_t$.
Suppose that the charging cost at time $t$ is a strictly convex increasing function of the total load, denoted by $f(s_t+l_t)$.
The convexity and increasing property of $f(s_t+l_t)$ reflects the fact that each unit of additional power demand becomes more expensive to obtain and make available to the consumer. For example, in the wholesale market, the instantaneous cost can be modeled as an increasing quadratic function of the instant load \cite{tang2014online,ma2013decentralized,he2012optimal}.
On the other hand, the convexity of $f(s_t+l_t)$ also captures the intent of reducing the load fluctuation over time \cite{gan2012optimal}.
Then the total cost over time $T$ is computed as $\sum_{t = 1}^{T} f(s_t+l_t)$.
In the ideal case, assume that $l_t, t_i^{(s)}, t_i^{(e)}$, and $d_i$ for all $t = 1, \cdots, T, i \in \mathcal{N}$ are known non-causally at the beginning of the system time. Then, the charging station can solve (\ref{model6}) and obtain the optimal charging rate, denoted by $x_{it}^*$ for all time $t$ and the optimal total cost, denoted by $\Psi_{1}$.
Such a solution is referred to as an ``optimal offline solution''.
\begin{subequations} \small
\label{model6}
\begin{align}
\Psi_{1} = & \min_{x_{it}} & & \quad \sum_{t = 1}^{T} f \left (\sum_{i \in \mathcal{I}(t)} x_{it} + l_t \right) \label{model7}
\\
& \text{s.t.} & & \quad \sum_{t = t_i^{(s)}}^{t_i^{(e)}} x_{it} = d_i, \forall i \in \mathcal{N}, \label{model8}
\\
&  & & \quad x_{it} \geq 0, \forall t = t_i^{(s)}, \cdots, t_i^{(e)}, \forall i \in \mathcal{N}. \label{model9}
\end{align}
\end{subequations}
In particular, the optimal total charging rate, denoted by $s_t^*$, is defined as $s_t^* = \sum_{i \in \mathcal{I}(t)} x_{it}^*$.
Note that there are in total $O(T |\mathcal{I}(t)| )$ variables in (\ref{model6}), where $|\mathcal{I}(t)|$ denotes the cardinality of the set $\mathcal{I}(t)$. This number can be quite large when the number of cars present at each time slot, $|\mathcal{I}(t)|$, is large.  Next, we propose an equivalent transformation of (\ref{model6}) that drastically reduces the number of variables. In particular, the following Theorem \ref{theorm1} shows that as long as we find the optimal $s_t^* ~\forall t$, the optimal $x_{it}^* ~\forall i, t$ can be obtained by earliest deadline first (EDF) scheduling.
\begin{theorem}\label{theorm1}
If a set of $s_t$'s satisfy the following inequality for all $n = 1, \cdots, T$
\begin{equation}\small\label{theorem:eq:s}
\sum_{t=1}^n \sum_{ i \in \{ i|t_i^{(e)} = t \} } d_i \leq \sum_{t=1}^{n} s_t  \leq \sum_{t=1}^n \sum_{ i \in \{ i|t_i^{(s)} = t \} } d_i,
\end{equation}
then there exists at least a set of $x_{it}$'s that is feasible to (\ref{model6}).
One such set of $x_{it}$'s can be obtained by EDF scheduling, which charges the PEV $i \in \mathcal{I}(t)$ with the earliest deadline at a rate $s_t$ at each time $t$.
Moreover, when $s_t=s_t^*$,
the set of $x_{it}$'s obtained by EDF scheduling are the optimal solution, $x_{it}^*$, to (\ref{model6}).
\end{theorem}
\emph{Proof:} Please see the detailed proof in Appendix \ref{appendix:proof:theorem1}.
To see Theorem \ref{theorm1}, note that \eqref{theorem:eq:s} implies that the total energy charged by any time slot $n$ is no less than the total charging demand that must be satisfied by time $n$ and no more then the total charging demand of PEVs which have arrived up to time $n$. On the other hand, by EDF scheduling, PEVs with earlier deadlines must be fully charged before those with later deadlines can be charged. Thus, \eqref{theorem:eq:s} guarantees the fulfillment of the charging demands of each individual PEV.
With Theorem \ref{theorm1}, we can transform (\ref{model6}) to the following equivalent problem with $T$ variables.
\begin{subequations}\small 
\label{ws_model}
\begin{align}
\Psi_{1} = & \min_{s_t} & & \quad \sum_{t = 1}^{T} f(s_t+l_t) \label{ws_modela}
\\
& \text{s.t.} & & \quad \sum_{t=1}^{n} s_t \geq \sum_{j=1}^n \sum_{ i \in \{ i|t_i^{(e)} = j \} } d_i, \forall n = 1, \cdots, T. \label{ws_modelb}
\\
& & & \quad  \sum_{t=1}^{n} s_t \leq \sum_{j=1}^n \sum_{ i \in \{ i|t_i^{(s)} = j \} } d_i, \forall n = 1, \cdots, T. \label{ws_modelc}
\end{align}
\end{subequations}
The optimal solution $s_t^*$ to (\ref{ws_model}) has an interesting feature: it does not change with the cost function $f(s_t+l_t)$, as long as $f$ is strictly convex. Moreover, $s_t^*$ also minimizes the variance of total load subjecting to \eqref{ws_modelb} and \eqref{ws_modelc}, where the variance of total load is defined as $\sum_{t =1}^T (s_t + l_t - \frac{ \sum_{t =1}^T s_t+l_t }{T})^2$ \cite{chen2014distributional,gan2013realtime}.
This is proved in  Theorem \ref{theorem:convex}.

\begin{theorem}\label{theorem:convex}
The optimal solution $s_t^*$ to  (\ref{ws_model}) does not change with the cost function $f(.)$, as long as $f(.)$ is strictly convex. Moreover, $s_t^*$ is a load flattening solution that minimizes the variance of total load.
\end{theorem}
\emph{Proof:} Please see the detailed proof in Appendix \ref{appendix:proof:theorem:convex}.

\remark{In practice, a capacity constraint on $s_t+l_t$ is present for each $t$ due to the hardware limitations and security concerns. The constraint is omitted in our formulation for the following reason. Theorem \ref{theorem:convex} indicates that the optimal solution $s_t^*$ to (\ref{ws_model}) minimizes the variance of total load. That is, any other scheduling solution would have a higher peak load, and therefore is more likely to violate the capacity constraint. In this sense, the optimal solution $s_t^*$ to (\ref{ws_model}) is ``capacity optimal'' in the sense that if the optimal solution to Problem \ref{ws_model} (or equivalently Problem \ref{model6}) violates the capacity constraint, then there does not exist any other scheduling solutions that satisfy the capacity constraint. }


\subsection{Online PEV Charging problem}\label{subsec:onlineproblem}
For the online PEV charging problem, the charging schedule only depends on the statistic information of future load demand, the current based load and the remaining charging demands and deadlines of the PEVs that have arrived so far.
In contrast to the ideal case in the last subsection, in practice the charging station only knows the remaining charging demands and departure deadlines of the PEVs that have already arrived, as well as the current base load.
The online charging scheduling algorithm computes the charging rate $s_k$ at each time slot $k$ based on the known causal information and the statistical information of the unknown future demands.
The charging rate $s_k$, once determined, cannot be changed in the future.
Specifically, the remaining charging demand of PEV $i$ at time $k$ is given by $\hat{d}_i^k =d_i -\sum_{t = t_i^{(s)}}^{k-1} x_{it}.$
Note that, $\hat{d}_i^k = d_i$ for all PEVs that have not yet arrived by time $k-1$.
A close look at (\ref{ws_model}) suggests that the charging schedule $s_t$ only depends on the total charging demand that needs to be finished before a certain time, but not the demand due to individual PEVs.
Thus, for notational simplicity, we define $\tilde{d}_{t}^k = \sum_{ i \in \{ i|t_i^{(e)}=t \} }\hat{d}_i^k, \forall t = k, \cdots, T,$
as the total unfinished charging demand at time $k$ that must be completed by time $t$.
With this, we define the state of system at time $t$ as
\begin{equation}\small
\label{def_state}
\mathbf{D}_t = [l_t, \tilde{d}_t^t, \tilde{d}_{t+1}^t, \cdots, \tilde{d}_T^t ],
\end{equation}
where $l_t$ is the base load at time $t$, $\tilde{d}_{t'}^t$ is the total unfinished charging demand at time $t$ that must be completed by time $t'$.
Let $\boldsymbol{\xi}_{t}$ represent the random arrival events at time $t$. $\boldsymbol{\xi}_t$ is defined as
\begin{equation}\small 
\label{def_xi}
\boldsymbol{\xi}_t = [\iota_t, \eta_{t}^{t}, \eta_{t+1}^{t}, \cdots, \eta_{e_t}^{t}],
\end{equation}
where $\iota_t$ is the base load at time $t$, $\eta_{t'}^{t}$ is the total charging demand that arrive at time $t$ and must be fulfilled by time $t'$, $e_t$ is the latest deadline among the PEVs that arrive at time $t$.
Then, the state transition, defined as
\begin{equation}\small
\mathbf{D}_{t+1} := g( s_t, \mathbf{D}_t, \boldsymbol{\xi}_{t+1} ),
\end{equation}
is calculated as follows: 
\begin{equation}\small
l_{t+1} = \iota_{t+1}
\end{equation}
and
\begin{equation}\small
\tilde{d}_{t'}^{t+1} = \left[\tilde{d}_{t'}^t -  \left[s_t - \sum_{j = t}^{t'-1} \tilde{d}_j^t \right]^+ \right]^+ + \eta_{t'}^{t+1}, \forall t' = t+1, \cdots, T.
\end{equation}
Here, $[x]^+ = \max \{ x, 0 \}$.
With the above definitions of system state and state transition, we are now ready to rewrite (\ref{ws_model}) into the following finite-horizon dynamic programming problem.
\begin{subequations}\small\label{msp_model1} 
\begin{align}
Q_k(\mathbf{D}_k)= &\min_{s_k} & &  f(s_k+l_k) + \mathbf{E}_{\boldsymbol{\xi}_{k+1}}[ Q_{k+1}( g( s_k, \mathbf{D}_k, \boldsymbol{\xi}_{k+1} ) ) ] \label{msp_model1a}
\\
&\text{s. t. }  & & \tilde{d}_k^k \leq s_k \leq \sum_{t = k}^T \tilde{d}_t^k, \label{msp_model1b}
\end{align}
\end{subequations}
where $Q_{k}( \mathbf{D}_{k} )$ is the optimal value of the dynamic programming at time $k$.
The left side of (\ref{msp_model1b}) ensures all charging demands to be satisfied before their deadlines. The right side of (\ref{msp_model1b}) implies that the total charging power up to a certain time cannot exceed the total demands that have arrived up to that time.
By slight abuse of notation, in the rest of the paper we denote the optimal solutions to both the online and offline problems as $s_k^*$, when no confusion arises. The actual meaning of $s_k^*$ will be clear from the context.
Suppose that $s_k^*$ is the optimal solution to (\ref{msp_model1}) at stage $k$.
Then, the total cost at the end of system time, denoted by $\Psi_{2}$, is provided by
\begin{equation}\small\label{eq:psi_on}
\Psi_{2} = \sum_{k = 1}^T f(s_k^*+l_k).
\end{equation}
Note that (\ref{msp_model1a}) comprises nested expectations with respect to the random PEV arrivals in the future time stages. Except for few special cases, it is hard to provide the closed-form of the optimal solution to (\ref{msp_model1}).
On the other hand, (\ref{msp_model1}) can be solved by the standard numerical methods, such as backward reduction and the sample average approximation (SAA) based on Monte Carlo sampling techniques \cite{shapiro2009lecture,birge1997introduction,defourny2011multistage,maggioni2012analyzing}.
These algorithms typically incur a computational complexity that grows exponentially with both the time span $T$ and the dimensions of state and decision spaces. Note that (\ref{msp_model1}) involves continuous state and decision spaces. Discretization of these spaces leads to a curse of dimensionality, rendering the computational complexity prohibitively high.

\section{MPC-based Online Charging Algorithm}\label{sec:online algorithm}
In view of the extremely high complexity of standard numerical methods, we are motivated to obtain a near-optimal solution by solving a much simpler problem, which replace all exogenous random variables by their expected values. This is referred to as the expected value problem \cite{birge1997introduction,defourny2011multistage,maggioni2012analyzing} or the MPC approach \cite{rao2014smart,chen2014distributional,gan2013realtime,bansal2014plug} in the literature.
Notice that the first-order moment, i.e., expectation, of a random process is much easier to estimate than the other statistics, e.g., variance or the probability distribution. Thus, the assumption of knowing the expected values is weaker than the assumptions in other EV-charging algorithms \cite{mean2008}, which assume that the probability distributions of the random processes are known.

Instead of solving the problem using generic convex optimization approaches, we propose a low-complexity online Expected Load Flattening (ELF) algorithm by exploring the load flattening feature of the optimal solution to the expected value problem, as shown in Section \ref{subsec:algorithm description}.
Section \ref{subsec:optimality} provides the theoretical analysis of the performance gap between the optimal solution to the expected value problem and the optimal solution to (\ref{msp_model1}).

\subsection{Algorithm Description}\label{subsec:algorithm description}
Denote the expectation of $\boldsymbol{\xi}_t$ as $\boldsymbol{\mu}_t= [\nu_t, \mu_t^t, \cdots, \mu_T^t],$
where $\nu_t = \mathbf{E}[\iota_t], \mu_{t'}^t = \mathbf{E}[\eta_{t'}^t], \forall t' = t, \cdots,T.$
Replacing $\boldsymbol{\xi}_t$ in (\ref{msp_model1}) with $\mu_t$, we obtain the following deterministic problem:
\begin{subequations}\small\label{dp_model5}
\begin{align}
& \min_{s_k} & &  f(s_k + l_k)+\sum_{t = k+1}^T f (s_t+\nu_t) \label{dp_model5a}
\\
&\text{s. t. }  & &  \sum_{t = k}^j s_t \geq \sum_{t = k}^{j}\tilde{d}_t^k + \sum_{m=k+1}^j \sum_{n = m}^j \mu_{n}^{m}, \forall j = k, \cdots, T, \label{dp_model5b}
\\
& & &  \sum_{t = k}^j s_t  \leq \sum_{t = k}^{T}\tilde{d}_t^k + \sum_{m=k+1}^j \sum_{n = m}^{e_m} \mu_{n}^{m}, \forall j = k, \cdots, T. \label{dp_model5c}
\end{align}
\end{subequations}
In each time $k$, we solve problem (\ref{dp_model5}) and obtain the optimal charging solution $s_k^*$.  Then, problem (\ref{dp_model5}) is resolved with the updated $\tilde{d}_t^{k}$ according to the realization of the PEVs arrived in next time. So on
and so forth, we obtain the optimal charging solution $s_k^*$ for time stage $k = 2, \cdots, T$. The total cost at the end of system time, denoted by $\Psi_{3}$, is defined as
\begin{equation}\small\label{eq:psi_po}
\Psi_{3} = \sum_{k = 1}^T f(s_k^* +l_k),
\end{equation}
where $s_k^*$ is the optimal solution to (\ref{dp_model5}) at time stage $k$.
The solution to (\ref{dp_model5}) is always feasible to (\ref{msp_model1}) in the sense that it always guarantees fulfilling the charging demand of the current parking PEVs before their departures.
This is because the constraints of $s_k$ in (\ref{msp_model1}) are included in (\ref{dp_model5}).

Due to the convexity of $f(\cdot)$, the optimal solution is the one that flattens the total load as much as possible.
By exploiting the load flattening feature of the solution, we present in Algorithm \ref{alg:online} the online ELF algorithm that solves (\ref{dp_model5}) with complexity $O(T^3)$.
The online ELF algorithm have a lower computational complexity than generic convex optimization algorithms, such as the interior point method, which has a complexity $O(T^{3.5})$ \cite{ye1997}.
Notice that similar algorithms have been proposed in the literature of speed scaling problems \cite{yao1995a,bansal2007speed} and PEV charging problems \cite{tang2014online}. The optimality and the complexity of the algorithm have been proved therein, and hence omitted here. The algorithm presented here, however, paves the way for further complexity reduction to $O(1)$ in Section \ref{sec:case_study}.
For notation brevity, we denote in the online ELF algorithm 
\begin{equation}\small \label{dp2}
\begin{aligned}
\bar{d}_{t''}^{t'}
= \begin{cases}
\tilde{d}_{t''}^{t'}, & \text{for $ t'' = k, \cdots, T, t' = k,$} \\
\mu_{t''}^{t'}, & \text{for $t'' = t', \cdots, T, t' = k+1, \cdots, T$}.
\end{cases}
\end{aligned}
\end{equation}
The key idea of online ELF algorithm is to balance the charging load among all time slots $k, \cdots, T$. 
Specifically, step 3 - 5 is to search the time interval $[i^*, j^*]$ that has the maximum load density among current time slots and record the maximum load density. The time slots with maximum load density are then deleted, and the process is repeated until the current time slot $k$ belongs to the maximum-density interval, i.e., $i^* = k$.
\begin{algorithm}\label{alg:online} 
\small
\SetAlgoLined
 \SetKwData{Left}{left}\SetKwData{This}{this}\SetKwData{Up}{up}
 \SetKwRepeat{doWhile}{do}{while}
 \SetKwFunction{Union}{Union}\SetKwFunction{FindCompress}{FindCompress}
 \SetKwInOut{Input}{input}\SetKwInOut{Output}{output}
 \Input{$\mathbf{D}_k, \boldsymbol{\mu}_t, t = k+1, \cdots, T$}
 \Output{$s_k$}
 initialization $i = 0, j = 0$\;
\Repeat{$i^* = k$}{
For all time slot $i = k, \cdots, T, j = i, \cdots, T$, compute
\begin{equation}\small \label{algorithm:eq}
i^*, j^* = \arg \max_{k \leq i \leq j \leq T} \{ \frac{\sum_{t'=i}^j (\sum_{t''=t'}^j \bar{d}_{t''}^{t'}+\nu_{t'})}{j-i+1} \}.
\end{equation}

Set
\begin{equation}\small
y^* = \frac{\sum_{t'=i^*}^{j^*} (\sum_{t''=t'}^{j^*} \bar{d}_{t''}^{t'}+\nu_{t'})}{j^*-i^*+1}.
\end{equation}

Delete time slot $i^*, \cdots, j^*$ and relabel the existing time slot $t > j^*$ as $t - j^* + i^* - 1$.
}
Set $s_k = y^* - l_k$.
\caption{Online ELF Algorithm} \label{onlinealgorithm}
\end{algorithm}

\subsection{Optimality Analysis}\label{subsec:optimality}
In this subsection, we analyze the optimality of the solution to (\ref{dp_model5}).
Notice that MPC approximates the non-causal random variables by their expected values regardless of their distribution functions. As a result, such approximation may lead to unacceptably large performance loss, depending on the distribution of the random variables. Therefore, the MPC approximation is not always justifiable.
A well-accepted metric, \emph{Value of the Stochastic Solution} (VSS) is adopted to evaluate optimality gap between the optimal online solution and the solution to the expected value problem \cite{birge1997introduction,defourny2011multistage,maggioni2012analyzing}.
Previous work, e.g., \cite{defourny2011multistage,maggioni2012analyzing}, mainly evaluates VSS using numerical simulations. Whereas in our analysis, we show that VSS is always bounded regardless of the distribution of the future EV charging demands. This provides a strong theoretical justification of adopting MPC approximation to solve our problem.

Let $\Xi$ denote a scenario, which is defined as a possible realization of the sequence of random load demand \cite{shapiro2009lecture},
\begin{equation}\small
\Xi = [\boldsymbol{\xi}_2, \boldsymbol{\xi}_3, \cdots, \boldsymbol{\xi}_T].
\end{equation}
Here, we treat $\boldsymbol{\xi}_1$ as deterministic information since the demand of PEVs arrived at the first stage is known by the scheduler.
Let $\Phi_{1}, \Phi_{2}$ and $\Phi_{3}$ be the expectation of the optimal value of the offline problem \eqref{ws_model}, the online problem (\ref{msp_model1}) and the expected value problem \eqref{dp_model5}, respectively, where the expectation is taken over the random scenarios. That is,
\begin{equation}\small
\Phi_{1} = \mathrm E_{\Xi} \left[ \Psi_{1} (\Xi) \right], \Phi_{2} = \mathrm E_{\Xi} \left[ \Psi_{2}(\Xi) \right], \Phi_{3} = \mathrm E_{\Xi} \left[ \Psi_{3}(\Xi) \right].
\end{equation}
It has been proved previously \cite{birge1997introduction,defourny2011multistage} that
\begin{equation}\small\label{eq:eoff:eonopt:eonpoly}
\Phi_{1} \leq \Phi_{2} \leq \Phi_{3}.
\end{equation}
To assess the benefit of knowing and using the distributions of the future outcomes, the VSS is defined as
\begin{equation}\small
\text{VSS} =\Phi_{3} - \Phi_{2}.
\end{equation}
To show that the online ELF algorithm yields a bounded VSS, we need to bound $\Phi_{3}$ and $\Phi_{2}$.
Generally, it is hard to calculate $\Phi_{2}$ or analyze the lower bound of $\Phi_{2}$ directly \cite{defourny2011multistage,maggioni2012analyzing}.
Thus, we choose to analyze the lower bound of $\Phi_{1}$ instead, since (\ref{eq:eoff:eonopt:eonpoly}) shows that the lower bound of $\Phi_{1}$ is also the bound of $\Phi_{2}$.
In what follows, we will show the lower bound of $\Phi_{1}$ in Proposition \ref{proposition:ev_lowerbound} and the upper bound of $\Phi_{3}$ in Proposition \ref{proposition:bound:onpoly}.

\begin{proposition}\label{proposition:ev_lowerbound}
\begin{equation}\small
\Phi_{1} \geq T f \left(\frac{\sum_{t = 1}^{e_1}\tilde{d}_t^1 +  \sum_{t = 2}^T \sum_{j = t}^{e_t} \mu_t^j + \sum_{t=1}^T \nu_t}{T} \right).
\end{equation}
\end{proposition}
\emph{Proof:}  Please see the detailed proof in Appendix \ref{appendix:proof:proposition:ev_lowerbound}.

Let $\mathcal{O}(t)$  be the set that $\mathcal{O}(t) = \{(m,n)|  e_m \geq t, m = 1, \cdots, t, n = t, \cdots, e_m \}.$
Then, we show that $\Phi_{3}$ is bounded by Proposition \ref{proposition:bound:onpoly}.
\begin{proposition}\label{proposition:bound:onpoly}
For any distribution of $\boldsymbol{\xi}_t, t = 1, \cdots, T$, there is
\begin{equation}\small
\Phi_{3} \leq \mathrm E \left [ \sum_{t=1}^T f \left( \sum_{(m,n) \in \mathcal{O}(t) } \eta_n^m + \iota_t \right) \right].
\end{equation}
\end{proposition}
\emph{Proof:}  Please see the detailed proof in Appendix \ref{appendix:proof:proposition:bound:onpoly}.

Now, we are ready to present Theorem \ref{theorem:bound:vss}, which states that the VSS is bounded for any distribution of random variables.
\begin{theorem}\label{theorem:bound:vss}
For any distribution of random vector $\boldsymbol{\xi}_t$, $t = 1, \cdots, T, n = t, \cdots, T,$
there is
\begin{equation}\small\label{eq:vss:bound}
\text{VSS} \leq \mathrm E \left [ \sum_{t=1}^T f \left (\sum_{(m,n) \in \mathcal{O}(t) } \eta_n^m + \iota_t \right ) \right] - T f \left ( \frac{\Gamma}{T} \right ),
\end{equation}
where $\Gamma=\sum_{t = 1}^{e_1}\tilde{d}_t^1 +  \sum_{t = 2}^T \sum_{j = t}^{e_t} \mu_t^j + \sum_{t=1}^T \nu_t$.
\end{theorem}

Theorem \ref{theorem:bound:vss} can be easily derived by Proposition \ref{proposition:ev_lowerbound} and Proposition \ref{proposition:bound:onpoly}.
In practice, the performance gap between the online ELF algorithm and the optimal online algorithm is often much smaller than the bound of VSS. This will be elaborated in the numerical results in Section \ref{sec:simulation}.





\section{Online ELF Algorithm under First-order Periodic Random Processes}\label{sec:case_study}
Notice that the complexity of $O(T^3)$ of online ELF algorithm mainly comes from step 3, which exhaustively searches the maximum-density period $[i^*, j^*]$ over all subintervals within $[k, T]$.
When the random arrival process is first-order periodic stochastically \footnote{The first-order periodic stochastic process is defined as a stochastic process whose first-order moment, i.e., mean is periodic. That is,  the mean of the random arrival events, i.e., $\boldsymbol{\mu}_{t}$ is periodic. However, the actual realizations of the arrival events $\boldsymbol{\xi}_{t}$ are uncertain and not periodic. },
we argue that the searching scope can be reduced to one period from the whole system time $T$.
Thus, the complexity of step 3 is limited by the length of a period instead of $T$.
As a result, the complexity of the algorithm reduces from $O(T^3)$ to $O(1)$, implying that it does not increase with the system time $T$, and thus the algorithm is perfectly scalable.
In practice, the arrival process of the charging demands are usually periodic stochastically.
For example, the arrival of charging demands at a particular location is statistically identical at the same time every day during weekdays (or weekends).
the National Household Travel Survey (NHTS) 2009 gathers information about daily travel patterns of different types of households in 2009, and shows that the daily travel statistics (e.g., Average Vehicle Trip
Length, Average Time Spent Driving, Person Trips, Person Miles of Travel) are very similar for each weekday or weekend, but different between weekday and weekend \cite{nhts2009}.
In Section \ref{subsec:periodic}, we investigate the case when the random load demand process is first-order periodic. In Section \ref{subsec:stationary}, we provide a closed-form solution to \eqref{dp_model5} for a special case when the load demand process is the first-order stationary.

\subsection{First-Order Periodic Process}\label{subsec:periodic}
In this subsection, we consider the case when the arrival process is first-order periodic.
Specifically, the first-order periodic process means that the first-order moment (i.e., mean) of the random process is periodic. That is, at current time stage $k$, for all $t = k+1, \cdots, T,$ we have
\begin{equation}
\boldsymbol{\mu}_{t} = \mathbf E [ \boldsymbol{\xi}_{t} ] = \boldsymbol{\mu}_{t+p},
\end{equation}
where $\boldsymbol{\xi}_{t}$ is the random arrival events at time $t$, $\boldsymbol{\mu}_{t}$ is the expectation of $\boldsymbol{\xi}_t$, and $p$ is the length of period.
Then, instead of considering $\boldsymbol{\mu}_t$ for $t=k+1, \cdots, T$, we only need to consider $\boldsymbol{\mu}_t$ for one period, i.e., for $t=k+1, k+p$:
\begin{equation}\small
\begin{aligned}
&\boldsymbol{\mu}_{k+1} = [\nu_{k+1}, \mu_{k+1}^{k+1}, \mu_{k+2}^{k+1}, \cdots, \mu_{k+e_1}^{k+1}, 0, \cdots, 0], \\
&\vdots\\
&\boldsymbol{\mu}_{k+p} = [\nu_{k+p}, \mu_{k+p}^{k+p}, \mu_{k+p+1}^{k+p}, \cdots, \mu_{k+e_p}^{k+p}, 0, \cdots, 0].
\end{aligned}
\end{equation}
Here, $e_{n} \leq T, n = 1, \cdots, p$ is the maximum parking time for PEVs arriving at time $k+n$.
Specially, we define $\hat{e}$ as $\hat{e} = \max \{e_{k+1}, e_{k+2}, \cdots, e_{k+p} \}.$
We decompose the search region $\{i, j| i = k ,\cdots, T, j = i, \cdots, T \}$ into three sub-regions, defined as $\Pi_1 = \{i, j| i =k, j = k, \cdots, k+\hat{e} \}$, $\Pi_2 = \{i, j| i = k, j = k+\hat{e}+1, \cdots, T \}$ and $\Pi_3 = \{i, j| i =k+1, \cdots, T, j = i, \cdots, T \}$, respectively.
We denote by $\hat{X}, \hat{Y}, \hat{Z}$ the maximum densities of region $\Pi_1, \Pi_2, \Pi_3$, respectively.
Indeed, the largest of $\hat{X}, \hat{Y}, \hat{Z}$ is the maximum density of the interval $[i^*, j^*] \subseteq [k, T]$ over all possible pairs $i,j \in \{ i = k ,\cdots, T, j = i, \cdots, T \}$.
Let $[\hat{i}_1, \hat{j}_1], [\hat{i}_2, \hat{j}_2], [\hat{i}_3, \hat{j}_3]$ be the intervals with the maximum density over region $\Pi_1, \Pi_2, \Pi_3$, respectively.
By definition, $\hat{i}_1 = \hat{i}_2 = k.$
Similar to the stationary case, $\hat{X}$ can be calculated by searching $\hat{j}_1$ over $\{k, \cdots, k+\hat{e}\}$. That is,
\begin{equation}\small\label{value:x}
\hat{X} = \max_{k \leq t \leq k+\hat{e}}  \frac{\sum_{n = k}^{t}(\tilde{d}_n^k + \nu_n)+ \sum_{n=k}^{t}\sum_{m=n}^t \mu_m^n }{n-k+1}.
\end{equation}
Moreover, Lemma \ref{lemma:periodic} shows that $\hat{Y}$ and $\hat{Z}$ can be calculated once the maximum density of interval $[k+1, k+ \hat{e}]$ has been obtained.
First, we introduce some definitions which help to show Lemma \ref{lemma:periodic}.
Let $\Pi_4$ be a subset of $\Pi_3$, where $\Pi_4$ is defines as $\Pi_4=\{i, j| i =k+1, \cdots, k+ \hat{e}, j = i, \cdots, k+ \hat{e}\}$, and $[\bar{i}, \bar{j}]$ be the interval with maximum density of region $\Pi_4$, i.e.,
$\bar{i}, \bar{j} = \underset{k+1 \leq i \leq j \leq k+ \hat{e} } {\mathrm{arg~max}} ~ \frac{\sum_{n=i}^{j} (\sum_{m=n}^{k+e_{n}} \mu_{m}^{n} + \nu_{n}) }{j-i+1}.$

\begin{lemma}\label{lemma:periodic}
The maximum densities of $\Pi_2$ and $\Pi_3$ are calculated by
\begin{equation}\small\label{value:hat:y:z}
\hat{Y} = \frac{ \sum_{n=k}^{\hat{j}_2} (\sum_{m=n}^{k+e_{n}} \mu_{m}^{n} + \nu_{n}) }{ \hat{j}_2-k +1 },
\hat{Z} = \frac{ \sum_{n=\hat{i}_3}^{\hat{j}_3} (\sum_{m=n}^{k+e_{n}} \mu_{m}^{n} + \nu_{n}) }{ \hat{j}_3-\hat{i}_3 +1 },
\end{equation}
respectively, where $\hat{i}_3 = \bar{i}$,
\begin{equation}\small \label{eq:hat:j:2}
\begin{aligned}
\hat{j}_2
= \begin{cases}
\max\{\bar{j}, k+\hat{e}+1\}, &   \text{if~} \bar{j} < \bar{i}+p, \\
\bar{j}+(r-1)p, & \text{otherwise.}
\end{cases}
\end{aligned}
\end{equation}
and
\begin{equation}\small \label{eq:hat:j:2}
\begin{aligned}
\hat{j}_3
= \begin{cases}
\bar{j}, &   \text{if~} \bar{j} < \bar{i}+p, \\
\bar{j}+(r-1)p, & \text{otherwise.}
\end{cases}
\end{aligned}
\end{equation}
\end{lemma}
\emph{Proof:}  Please see the detailed proof in Appendix \ref{appendix:proof:lemma:periodic}.

Based on Lemma \ref{lemma:periodic}, we can modified the searching region of step $3$ of online ELF algorithm as follows:
\begin{itemize}
\item
if $\bar{j} < \bar{i}+p$, the interval with the maximum density during time stages $[k+1, T]$ is $[\bar{i}, \bar{j}]$. Then, in step $3$ of the online ELF algorithm, the search region of $i,j$ is reduced from $\{i,j| i = k, \cdots, T, j = i,\cdots, T \}$ to $\{i,j| i = k, \cdots, \bar{i}, j = i,\cdots, \bar{i}, \bar{j} \}$.
\item
If $\bar{j} \geq \bar{i}+p$, the interval with the maximum density during time stages $[k+1, T]$ is $[\bar{i}, \bar{j}+(r-1)p]$. Then, in step $3$ of the online ELF algorithm, the search region of $i,j$ can be reduced from $\{i,j| i = k, \cdots, T, j = i,\cdots, T \}$ to $\{i,j| i = k, \cdots, \bar{i}, j = i,\cdots, \bar{i}, \bar{j}+(r-1)p \}$.
\end{itemize}
As a result, the searching region of the online ELF algorithm is only related to $[k+1, k+\hat{e}]$ instead of $T$. Thus, the computational complexity of the online ELF algorithm is $O(1)$ instead of $O(T^3)$ under first-order periodic process.

\subsection{First-order Stationary Process}\label{subsec:stationary}
In this subsection, we show that the optimal solution to (\ref{dp_model5}) can be calculated in closed form if the arrival process is first-order stationary.
Due to the page limit, we only provide the main results here.
By first-order stationary, we mean that the statistical mean of $\boldsymbol{\xi}_t$, i.e., $\nu_t$ and $\mu_{t'}^t, t' = t, \cdots, T$ only depends on the relative time difference $\tau = t' - t$, but not the absolute value of $t$. We can then replace $\nu_t$ by $\nu$ and replace $\mu_{t'}^t$ by $\mu_\tau$, where $\tau = t'-t$.
Then, $\boldsymbol{\mu}_t$  is no longer a function of $t$, and can be represented as
\begin{equation}\small
\boldsymbol{\mu} = [\nu, \mu_1, \mu_2, \cdots, \mu_{\bar{e}}, 0, \cdots, 0],
\end{equation}
where $\bar{e}$ is the maximum parking time of a PEV.
We denote by $X, Y, Z$ the maximum densities of region $\Pi_1, \Pi_2, \Pi_3$, respectively.
Then, $X$ is calculated by
\begin{equation}\small\label{value:x}
\begin{aligned}
&X = \\
&\max_{k \leq n \leq k+\bar{e}} \left \{ \frac{\sum_{t = k}^{n}\tilde{d}_t^k + \sum_{j=1}^{n}(n-k-j+1)\mu_j + l_k - \nu }{n-k+1} + \nu\right \},
\end{aligned}
\end{equation}
and $Y, Z$ are provided in Lemma \ref{lemma:stationary}.
\begin{lemma}\label{lemma:stationary}
The maximum densities of $\Pi_2$ and $\Pi_3$ are achieved by setting $i_2 = k, j_2 = T, i_3 = k+1, j_3 = T$, and calculated by
\begin{subequations}\small\label{value:y:z}
\begin{align}
&Y = \frac{ \sum_{t = k}^{k+\bar{e}}\tilde{d}_t^k + \sum_{j=1}^{k+\bar{e}}(T-k-j+1)\mu_j + l_k -\nu }{ T-k+1 } + \nu,\\
&Z = \frac{ \sum_{j=1}^{k+\bar{e}}(T-k-j+1)\mu_j }{T-k} + \nu.
\end{align}
\end{subequations}
\end{lemma}
\emph{Proof:}  Please see the detailed proof in Appendix \ref{appendix:proof:lemma:stationary}.
The largest of $X, Y,$ and $Z$
is the maximum density of the interval $[i^*, j^*] \subseteq [k, T]$ over all possible pairs $i,j \in \{ i = k ,\cdots, T, j = i, \cdots, T \}$.
Specially, if $X$ or $Y$ is the largest one, then
$k$ is already contained in the maximum-density interval, and thus $X$ or $Y$ is the optimal charging rate  at time $k$.
On the other hand, if $Z$ is the largest, then the maximum-density interval, i.e., $[k+1, T]$, does not include $k$.
Following Algorithm~\ref{alg:online}, we will delete the maximum-density interval and repeat the process.
Now, time slot $k$ is the only remaining time slot after deletion.
This implies that all charging demands that have arrived by time slot $k$ should be fulfilled during time slot $k$.
These arguments are summarized in Proposition \ref{proposition:stationary}, which provides the closed form solution to \eqref{dp_model5}.
\begin{proposition}\label{proposition:stationary}
When the random load demand process is first-order stationary, the optimal charging schedule to (\ref{dp_model5}) is given by  the following close-form:
\begin{numcases}{s_k^* =} \label{solution:stationary}
X - l_k, &  ~if~ $ X = \max \{ X, Y, Z \} $, \label{sol:x}\\
Y - l_k,  & ~if~ $ Y = \max \{ X, Y, Z \}$, \label{sol:y}\\
\sum_{t = k}^{k+\bar{e}}\tilde{d}_t^k, & ~otherwise. \label{sol:z}
\end{numcases}
\end{proposition}

\section{Simulations}\label{sec:simulation}
In this section, we investigate the performance of the proposed online ELF algorithm through numerical simulations.
All the computations are solved in MATLAB on a computer with an Intel Core i3-2120 $3.30 GHz$ CPU and 8 GB of memory.
For comparison, we also plot the optimal solution to \eqref{msp_model1}, which is obtained by SAA method \cite{shapiro2009lecture}, and a heuristic solution by the online AVG algorithm \cite{he2012optimal}, which is obtained by charging each PEV at a fixed rate, i.e., its charging demand divided by its parking time.
Let the expected cost of AVG algorithm is denoted by $\Phi_{4}$.
Define the relative performance loss of ELF and AVG compared with the optimal online solution as
$\frac{\Phi_{3} - \Phi_{2}}{\Phi_{2}}$ and $\frac{\Phi_{4} - \Phi_{2}}{\Phi_{2}}$, respectively.
Similar to \cite{tang2014online} \cite{he2012optimal}, we adopt an increasing quadratic cost function in the simulations, i.e., $f(s_t+l_t) = (s_t+l_t)^2.$ Note that the cost function is increasing and strictly convex, since the load $s_t+l_t$ is always non-negative.

\subsection{Average Performance Evaluation}\label{subsec:performance}

\begin{figure}
        \centering
        \includegraphics[width=2.8in]{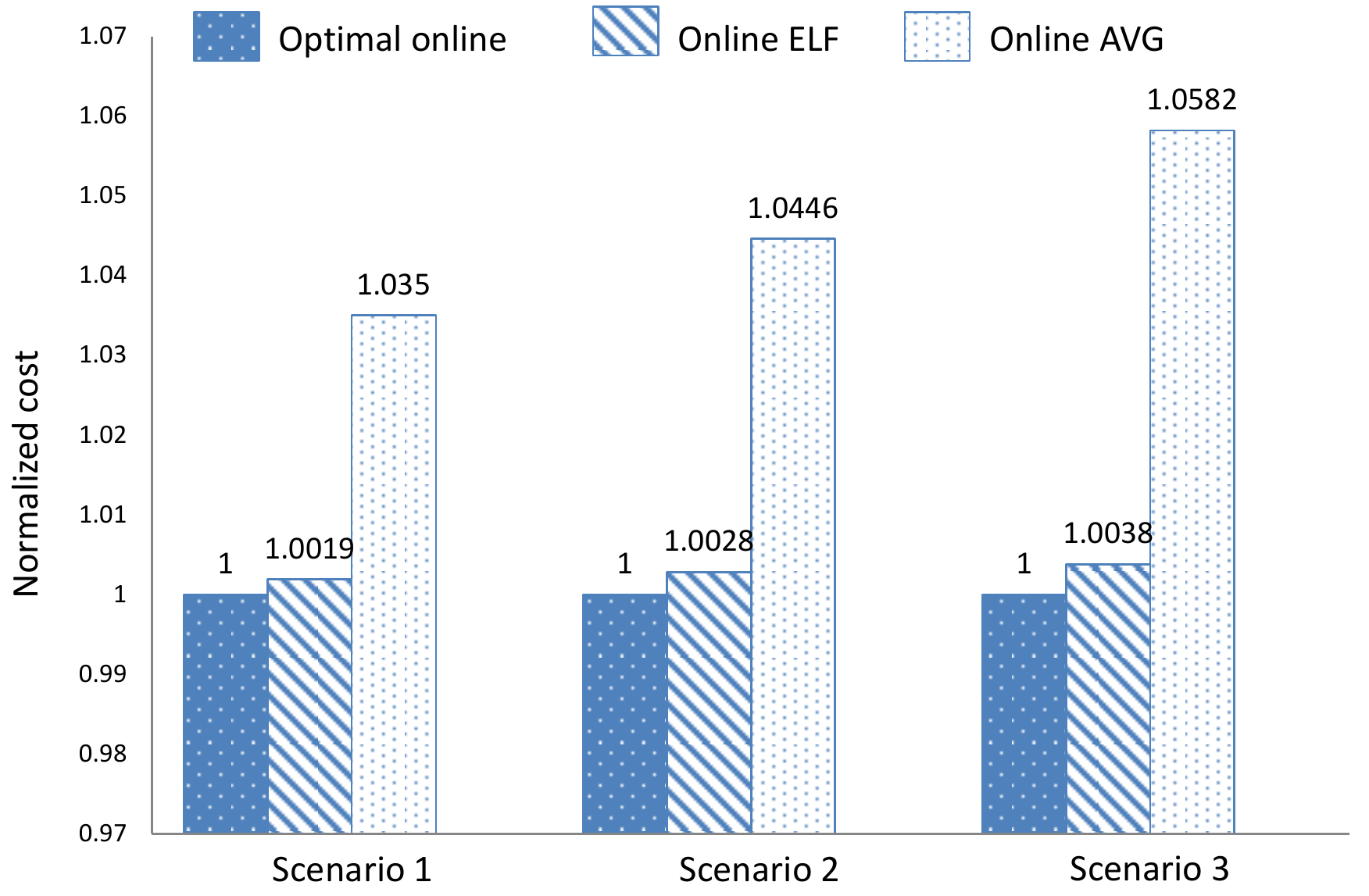}
        \caption{Normalized  costs of three algorithms in three scenarios.}
        \vspace{-0.7cm}
        \label{fig:normalized_cost}
\end{figure}

\begin{figure*}
\centering
\subfigure[Scenario 1: light traffic]{
    \includegraphics[width=0.31\textwidth]{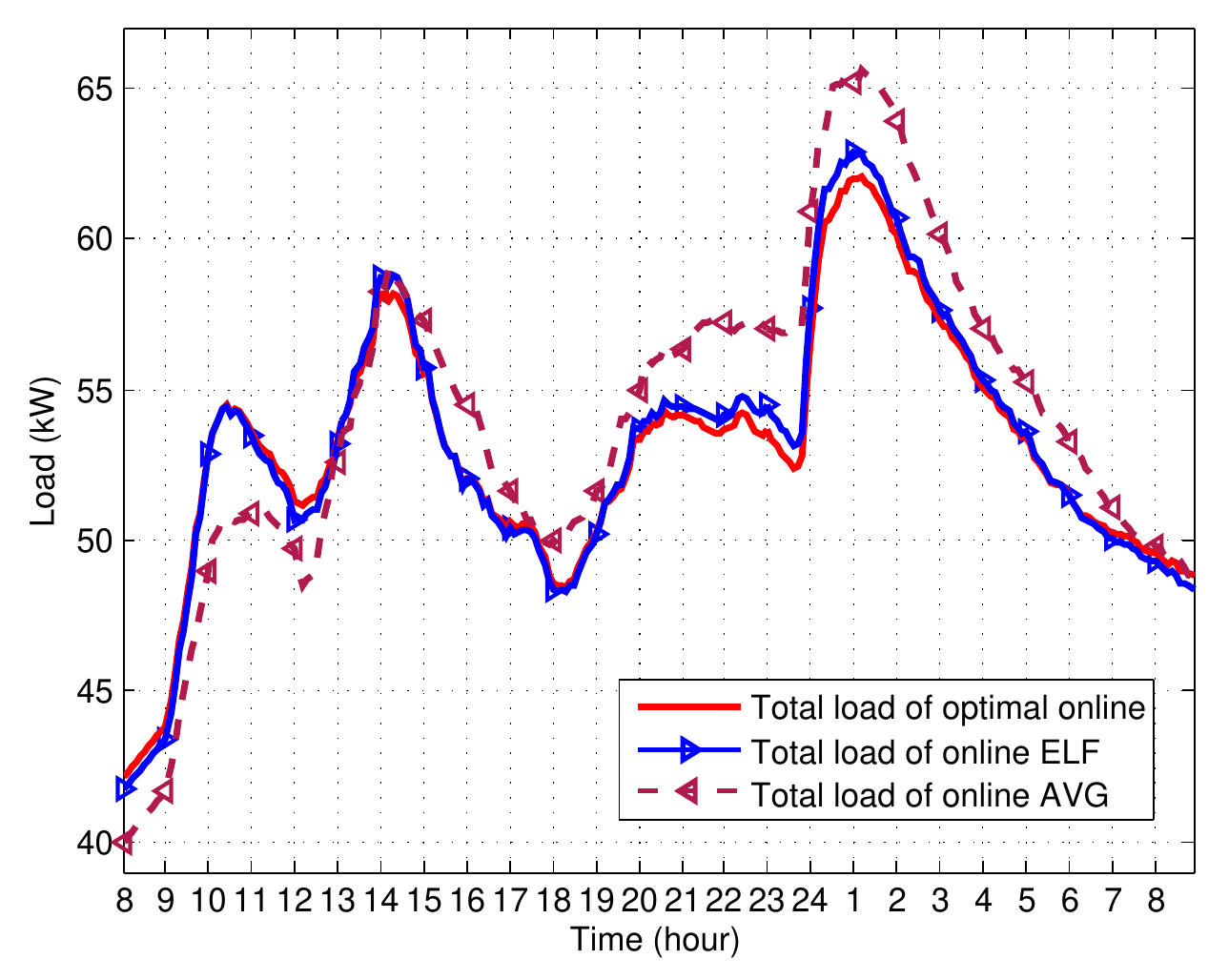}
        \label{fig:oneday_flat}
}
\subfigure[Scenario 2: moderate traffic]{
    \includegraphics[width=0.31\textwidth]{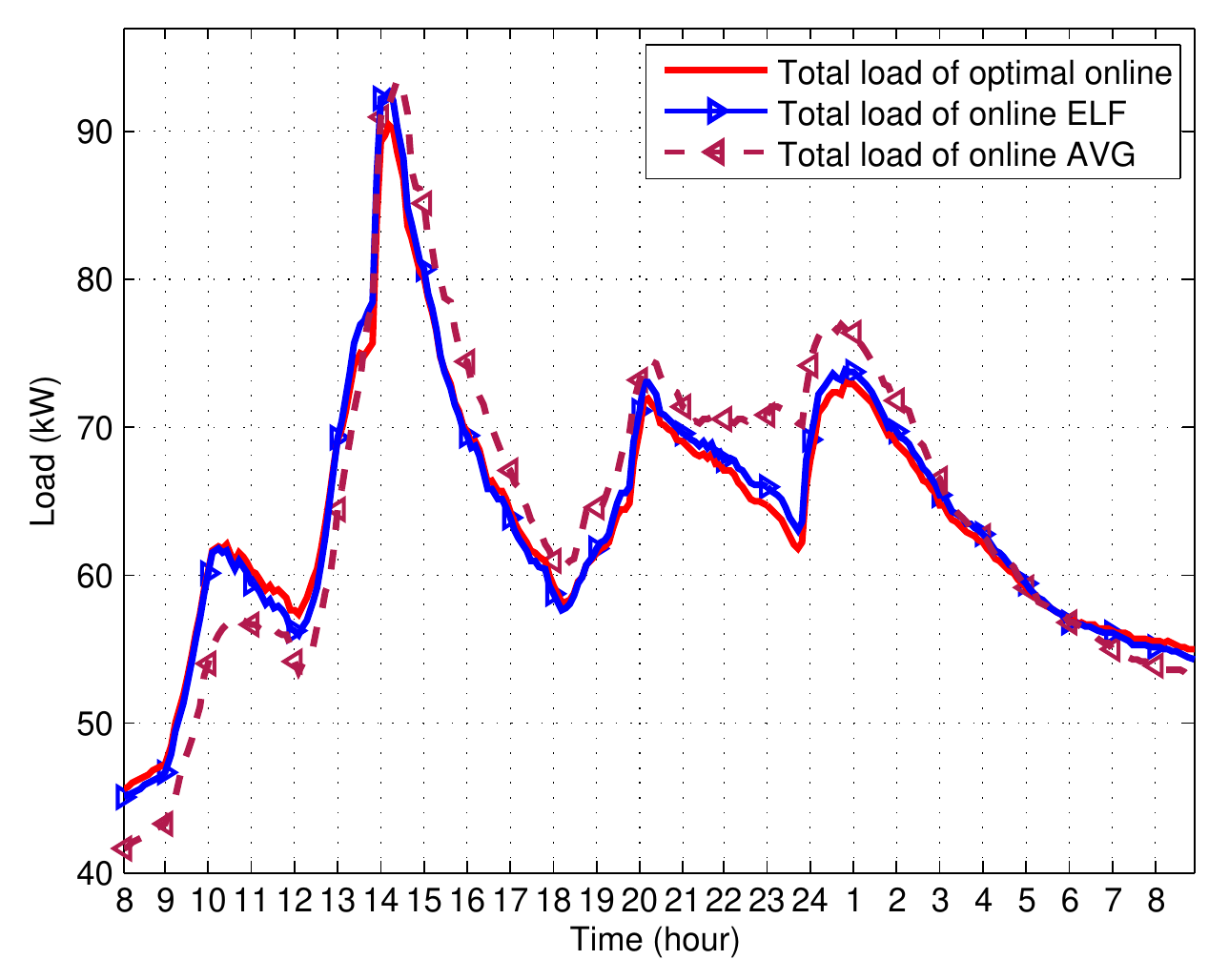}
        \label{fig:oneday_peak}
}
\subfigure[Scenario 3: heavy traffic]{
    \includegraphics[width=0.31\textwidth]{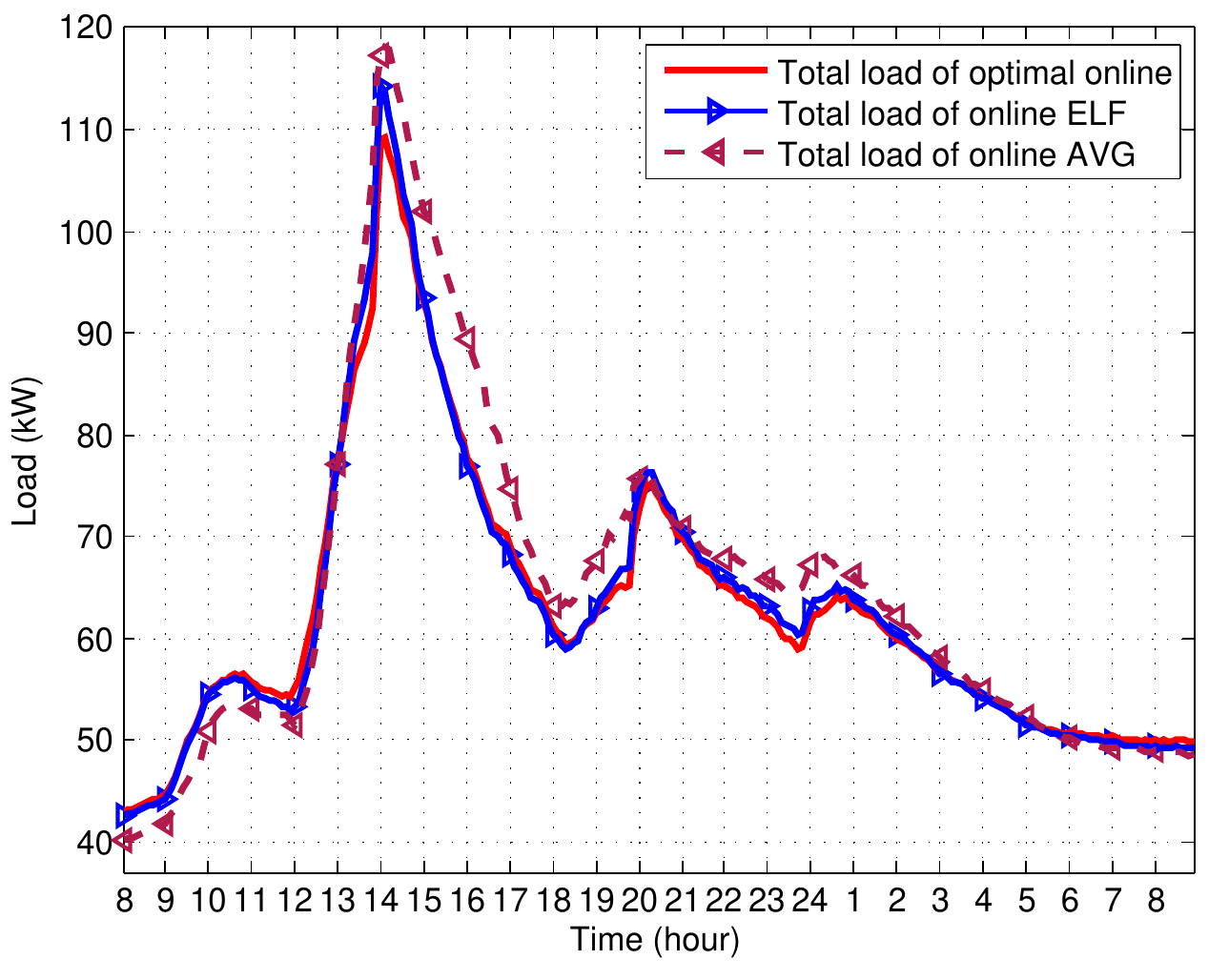}
        \label{fig:oneday_crowd}
}
\caption{
Base load and total load of three algorithms in three scenarios.
\label{fig:total_load}
}
\vspace{-0.3cm}
\end{figure*}

In this subsection, we evaluate the average performance of the online ELF algorithm under three different traffic patterns, i.e., light, moderate, and heavy traffics.
In particular, the system time is set to be $24$ hours, and each time slot lasts 10 minutes.
The PEV arrivals follow a Poisson distribution and the parking time of each PEV follows an exponential distribution \cite{alizadeh2014a,zhang2015an,chen2012iems}.
The mean arrival and parking durations of the three traffic patterns are listed in Table \ref{table:parameter2}.
The main difference lies in the arrival rates at the two peak hours, i.e. $12:00$ to $14:00$ and $18:00$ to $20:00$.
The settings of the peak hour match with the realistic vehicle trips in National Household Travel Survey (NHTS) 2009 \cite{nhts2009}.
Specially, the average number of total PEVs simulated in scenario 1, 2 and 3 are $104, 204$ and $304$, respectively.
We choose the base load profile of one day in the service area of South California Edison from \cite{gan2012optimal}.
Each PEV's charging demand is uniformly chosen from $[25, 35]kWh$.
Each point in Fig.~\ref{fig:normalized_cost}, Fig.~\ref{fig:total_load} and Table \ref{table:ratio1} is an average of $10^5$ independent instances of the scenarios listed in Table \ref{table:parameter2}.
In Fig.~\ref{fig:total_load}, the total loads $s_k^*$ are plotted over time.
We notice that the total load of the online ELF algorithm follows closely to that of optimal online algorithm, whereas that of the AVG algorithm has a larger gap from the optimal online algorithm.
The average costs normalized by that of the optimal online algorithm are plotted in Fig.~\ref{fig:normalized_cost}. Moreover, the VSS and the relative performance loss are listed in Table \ref{table:ratio1}. Both the figure and the table show that ELF performs very close to the optimal online algorithm.
The VSS and the relative performance loss are no more than $0.1536$ and $0.38\%$ respectively.
In contrast, the relative performance loss of AVG algorithm is up to $5.82\%$, which is more than 15 times of that of ELF algorithm.
The relative performance loss of the approximate online algorithm reflects the percentage of the extra cost compared with the optimal online algorithm.
Obviously, the performance loss is always the smaller the better.
For example, from the report of Rocky Mountain Institute, the average electricity rate is $11.2 cents/kWh$, and the average load for a charging station is $100 kW$ \cite{cost2009}. Then, the expected electricity cost of a charging station for one years is $\$ 967680$.  $6\%$ relative performance loss means AVG algorithm leads to $\$ 58060 $ extra cost, while the proposed online algorithm with $0.38\%$ relative performance loss leads to $\$ 3677 $ extra cost, or an annual saving of $\$ 54383$ compared with the AVG algorithm.

\begin{table}
\small		
        \caption{Parameter settings of the PEV traffic patterns}
        \centering
        \begin{tabular}{c c c c c}
        \toprule
         \multicolumn{1}{c}{\multirow {2}{*}{Time of Day}} & \multicolumn{3}{c}{Arrival Rate (PEVs/hour)} & \multicolumn{1}{c}{Mean Parking} \\
         \multicolumn{1}{c}{} & \multicolumn{1}{c}{S. 1} & \multicolumn{1}{c}{S. 2} & \multicolumn{1}{c}{S. 3} & \multicolumn{1}{c}{Time (hour)} \\
        \midrule
         $08:00 - 10:00$ &  $7$ & $7$ & $7$ & $10$ \\
         $10:00 - 12:00$ &  $5$ & $5$ & $5$ & $1/2$ \\
         $12:00 - 14:00$ &  $10$ & $35$ & $60$ & $2$ \\
         $14:00 - 18:00$ &  $5$ & $5$ & $5$ & $1/2$ \\
         $18:00 - 20:00$ &  $10$ & $35$ & $60$ & $2$ \\
         $20:00 - 24:00$ &  $5$ & $5$ & $5$ & $10$ \\
         $24:00 - 08:00$ &  $0$ & $0$ & $0$ & $0$ \\
         \bottomrule
        \end{tabular}
        \label{table:parameter2}
        \vspace{-0.4cm}
\end{table}

\begin{table}
\small
        \caption{Average performance comparison under three traffic patterns}
        \centering
        \begin{tabular}{c c c c c c }
        \toprule
        Scenario &  VSS  & $\frac{\Phi_{3} - \Phi_{2}}{\Phi_{2}}$ & $\frac{\Phi_{4} - \Phi_{2}}{\Phi_{2}}$  \\
        \midrule
        1 & 0.1178 & 0.19$\%$ & 3.50$\%$ \\
        2 & 0.1319 & 0.28$\%$ & 4.46$\%$ \\
        3 & 0.1536 & 0.38$\%$ & 5.82$\%$ \\
         \bottomrule
        \end{tabular}
        \label{table:ratio1}
        \vspace{-0.5cm}
\end{table}

\subsection{Complexity of The Online Algorithm ELF}

\begin{figure}
        \centering
        \includegraphics[width=2.7in]{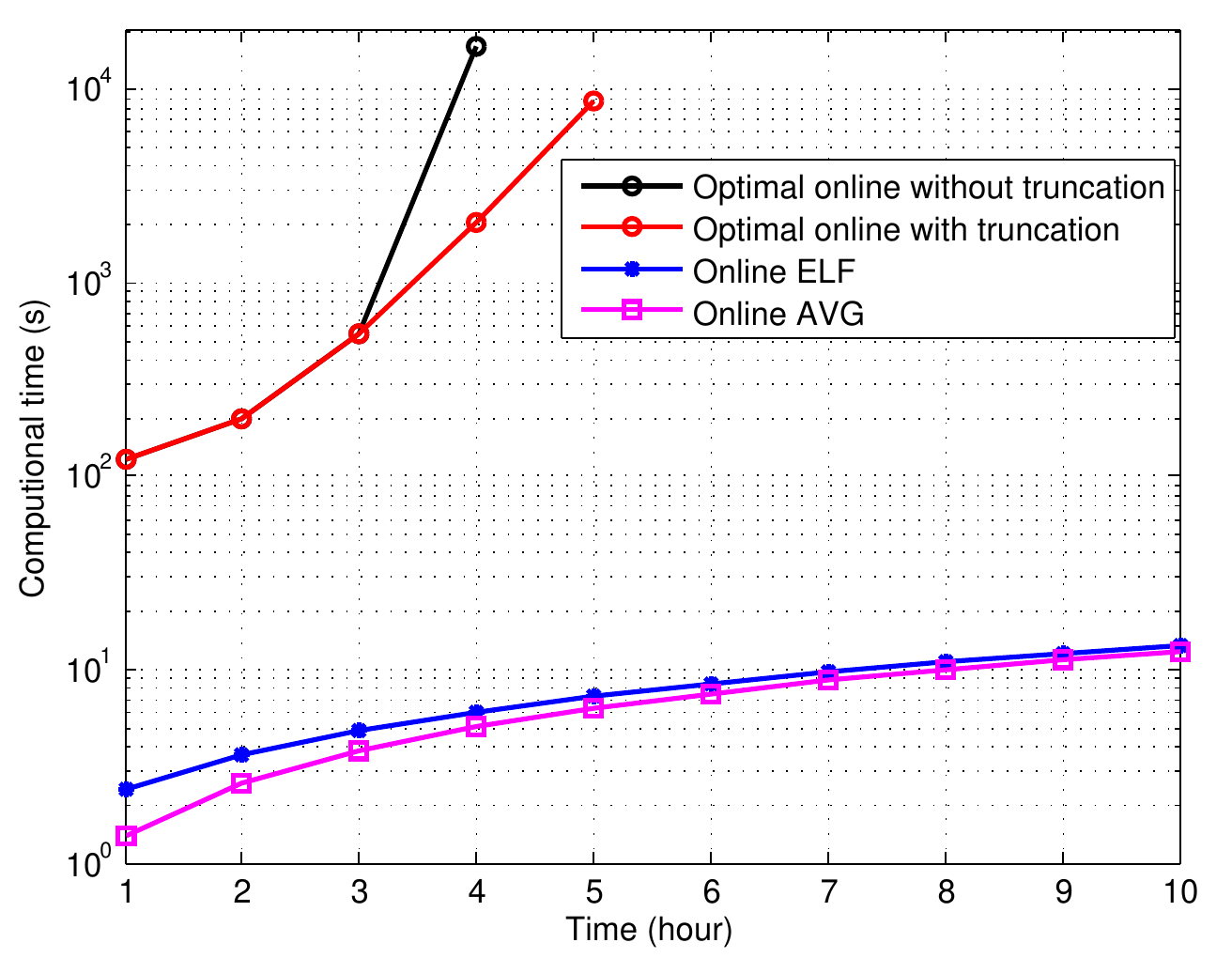}
        \caption{CPU computational time over the system time.}
        \vspace{-0.5cm}
        \label{fig:onlinecomplexity2}
\end{figure}

In this subsection, we verify the computational complexity of the online ELF algorithm and also compare the complexity of online ELF algorithm with that of the optimal online algorithm and online AVG algorithm, respectively.
We still adopt the SAA method as the optimal online algorithm.
Since the complexity of SAA is too high, truncation is often adopted in the SAA to reduce the complexity at a cost of performance loss \cite{birge1997introduction}. As such, we also simulate the complexity of truncated SAA with a truncation period of 3 hours.
Each PEV's charging demand is uniformly chosen from $[25, 35]kWh$, the arrivals follow a Poisson distribution and the parking time of each PEV follows an exponential distribution, where the arrival rate and mean parking durations are the same as that in the peak hours of scenario 2 in Section V.A.
We simulate $10$ cases in total, where the system time are set to be $1, 2, \cdots, 10$ hours.
For each algorithm, we record the CPU computational time  at each time stage and calculate the average CPU computational time as the sum of CPU computational times of all time stages divided by the number of time stages.
Each point in Fig. \ref{fig:onlinecomplexity2} is an average of $100$ independent instances.
Fig. \ref{fig:onlinecomplexity2} shows that the CPU computational time of both the online algorithm ELF and AVG almost grow linearly with the system time.
In contrast, for the SAA based online algorithm with or without truncation, the average CPU computational times grow very quickly as system time increases.
We notice that when system time increases to $4$ hours, the optimal online algorithm without truncation consumes more than $2$ hours and the optimal online algorithm with truncation consumes more than $30$ minutes. Meanwhile the proposed online algorithm ELF only takes several seconds.
It is foreseeable that the computational complexity of optimal online algorithm will become extremely expensive as we further increase the system time.

\subsection{Performance Comparison with Online Algorithm ORCHARD \cite{tang2014online}}
\begin{figure}
        \centering
        \includegraphics[width=2.7in]{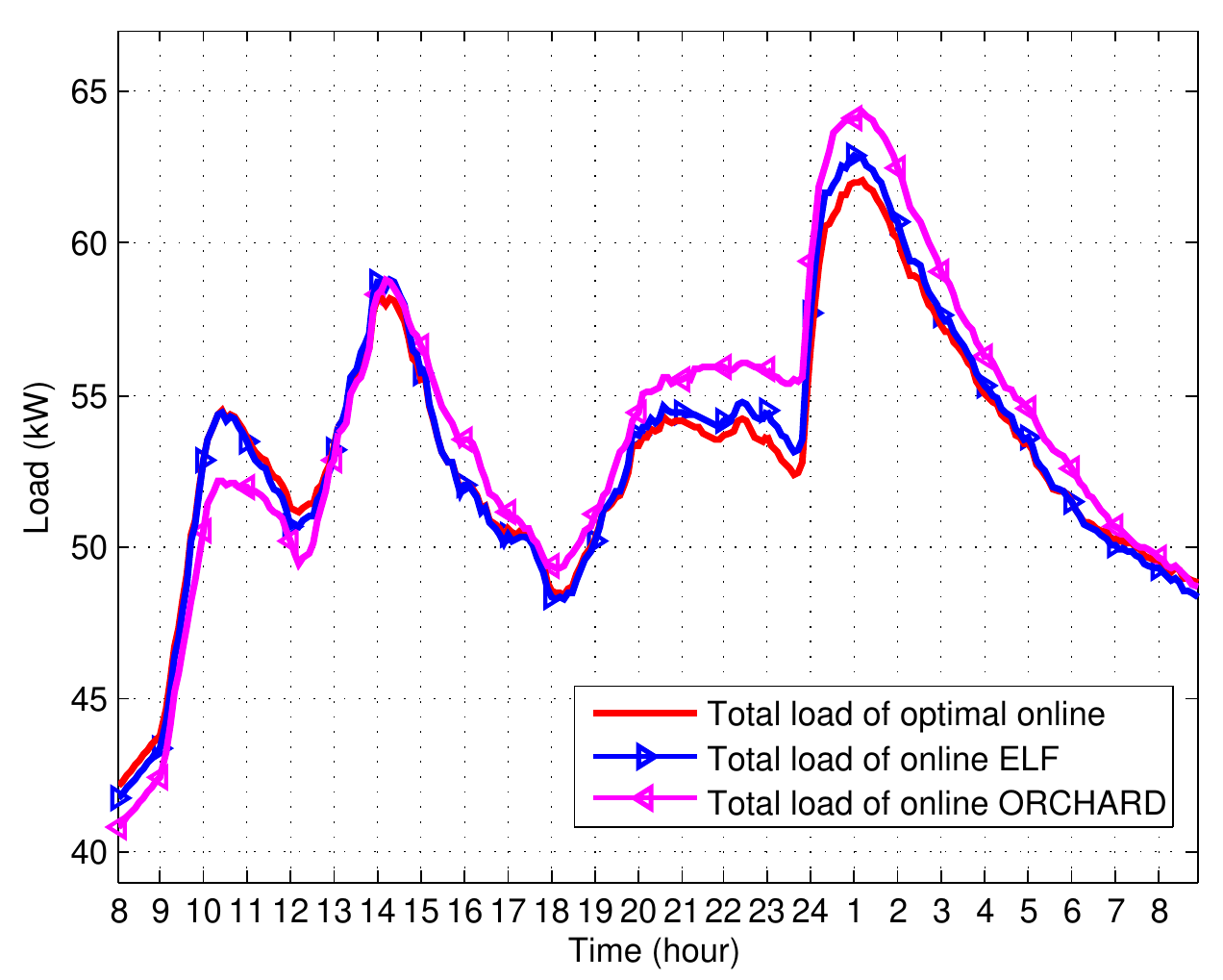}
        \caption{Load comparison of ORCHARD and ELF}
        \label{fig:loadcompare:orchard}
        \vspace{-0.5cm}
\end{figure}

\begin{figure}
        \centering
        \includegraphics[width=2.7in]{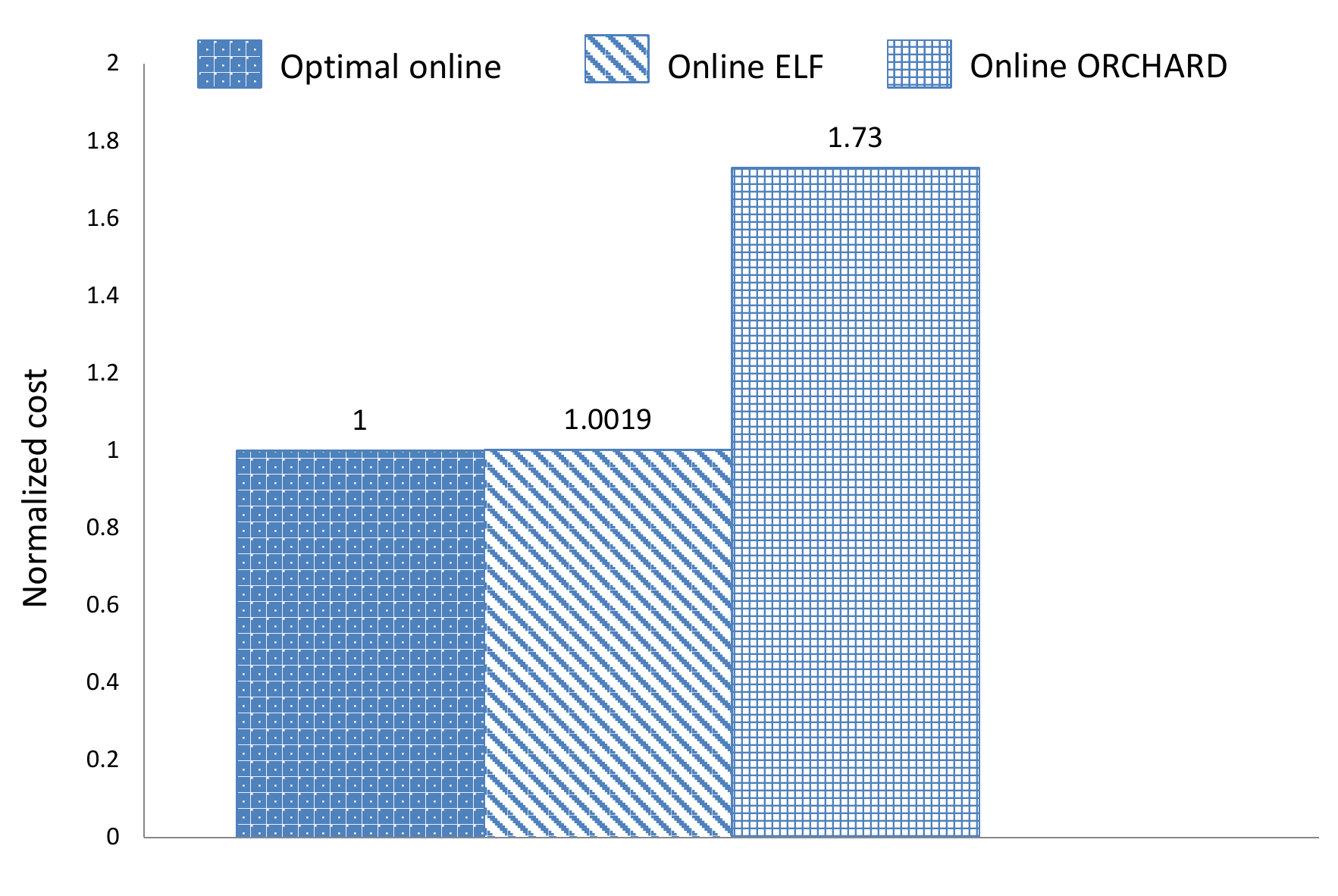}
        \caption{ Average performance ratios of online algorithm ORCHARD and ELF}
        \label{fig:ratiocompare:orchard}
        \vspace{-0.5cm}
\end{figure}

In this section, we compare the proposed online algorithm ELF with the online algorithm ORCHARD proposed in reference \cite{tang2014online} on the properties of the optimality and computational complexity.
First, we evaluate the average performance of the proposed online algorithm and online algorithm ORCHARD.
To facilitate the comparison, we also adopt the average performance of optimal online algorithm as a bench mark.
   For the system parameters, we use the default settings of scenario 1 in Section V.A.
    We simulate $10^5$ cases and plot the total load (the sum of the base load and the PEV charging load) over time for the three online algorithms in Fig. \ref{fig:loadcompare:orchard}.
    In addition, the average performance ratios normalized against the optimal online solution are shown in Fig.~ \ref{fig:ratiocompare:orchard} respectively.
    Fig.~\ref{fig:loadcompare:orchard} shows that compared with online algorithm ELF, the online algorithm ORCHARD always produces a larger gap from the optimal online solution. From Fig.~ \ref{fig:ratiocompare:orchard}, we can see that the proposed online algorithm ELF achieves a much lower expected cost than the online algorithm ORCHARD, which indicates that online algorithm ELF owns a better average performance than online algorithm ORCHARD.

\begin{figure}
        \centering
        \includegraphics[width=2.4in]{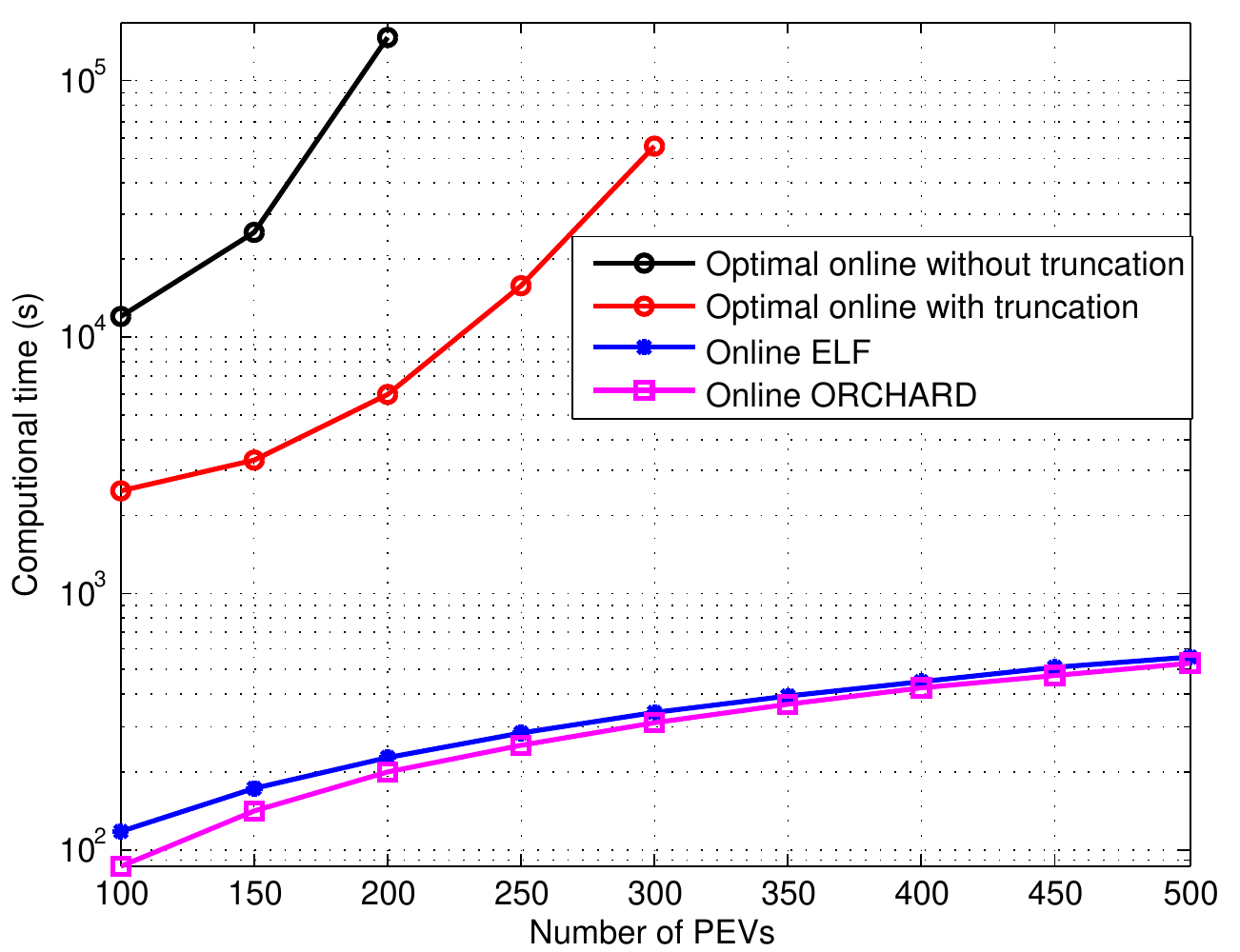}
        \caption{CPU computational time over the number of PEVs.}
        \vspace{-0.5cm}
        \label{fig:onlinecomplexity1}
\end{figure}

     To compare the computation complexity of online algorithm ORCHARD and ELF, we adopt the similar case study in Section V.D of reference \cite{tang2014online}.
     Specifically, we simulate the CPU computational time of online algorithms by varying the arrival rates of the PEVs during one day.
For the system parameter, we use the same settings as scenario 1 in Section V.A except the arrival rates, which are assumed to be the same during $8:00 - 18:00$ and $0$ after $18:00$. We vary
the arrival rate in $8 : 00 - 18 : 00$ from $10$ to $50$ (PEVs/hour) that leads to the average number of total PEVs during one day varies from $100$ to $500$.
For each specified average number of PEVs, we simulate the average performance of $10^8$ independent instances for the online ELF algorithm, the optimal online algorithm and the online ORCHARD algorithm, respectively, and record the average CPU computational times of each cases for the three algorithms, respectively. The results are plotted in Fig. \ref{fig:onlinecomplexity1}.
Fig. \ref{fig:onlinecomplexity1} shows that the average CPU computational time of the optimal online algorithm grows quickly as the number of total PEVs increases, while the average CPU computational time of the online ELF algorithm and the online ORCHARD algorithm grows slowly as the number of total PEVs increases. When the number of PEVs is $200$, the average CPU computational time of the optimal online algorithm without truncation is more than $24$ hours. Even for the the optimal online algorithm with truncation, the average CPU computational time is about $100$ minutes.
Whereas, the proposed online algorithm ELF only takes about $4$ minutes.
Fig. \ref{fig:onlinecomplexity1} also indicates that the computational time of online ELF and online ORCHARD as the number of the PEVs increases.

As a conclusion, the case study shows that the online algorithm ELF has a better average performance than online algorithm ORCHARD. It also indicates that the CPU computational time of online ELF and online ORCHARD are similar.

\section{Conclusions}\label{sec:conclusions}
In this paper, we formulate the optimal PEV charging scheduling problem as a finite-horizon dynamic programming problem.
Instead of adopting the standard numerical methods with high complexity, we provide a MPC-based online algorithm with $O(T^3)$-complexity.
We rigorously analyze the performance gap between the solution of the MPC-based approach and the optimal solution for any distribution of exogenous random variables.
Moreover, we show that the proposed algorithm can be made scalable under the first-order periodic process of load demand.
Besides, our analyses are validated through extensive simulations.

\appendix

\subsection{Proof of Theorem \ref{theorm1}:}\label{appendix:proof:theorem1}
We use the inductive method to show that through EDF scheduling, all the PEVs can be fulfilled charging before deadlines. For $n=1$, (\ref{theorem:eq:s}) becomes
\begin{equation}
 \sum_{ i \in \{ i|t_i^{(s)} = 1 \} } d_i \geq s_1 \geq \sum_{ i \in \{ i|t_i^{(e)} = 1 \} } d_i.
\end{equation}
Thus, by EDF scheduling, we can first satisfy the demand of PEVs whose deadline at time stage $1$. That is, for any PEV $i \in \{ i|t_i^{(e)} = 1 \}$, we set
\begin{equation}
x_{i1} = d_i.
\end{equation}
Assuming that for all time stage $m$, EDF scheduling can fulfill charge all the PEVs which depart at or before time stage $m$, i.e., there exists  at least a set of $x_{it}$'s that satisfy
\begin{subequations}\label{theorem:eq:s2}
\begin{align}
&\sum_{t = t_i^{(s)}}^{t_i^{(e)}} x_{it} = d_i,  \forall i \in \{ i|t_i^{(e)} \leq m \},\\
&x_{it} \geq 0, \forall t = t_i^{(s)}, \cdots, t_i^{(e)},  \forall i \in \{ i|t_i^{(e)} \leq m \}.
\end{align}
\end{subequations}
Since
\begin{equation}\small
\sum_{t=1}^{m} s_t \geq \sum_{t=1}^{m} \sum_{ i \in \{ i|t_i^{(e)} = t \} } d_i,
\end{equation}
then, $\sum_{t=1}^{m} s_t -  \sum_{t=1}^{m} \sum_{ i \in \{ i|t_i^{(e)} = t \} } d_i$ represents the amount of power which is outputted from the charging station during time stage $1, \cdots, m$ and charged to the PEVs with deadline after time stage $m$.
By EDF scheduling, once the PEVs which depart at time $m$ have been fulfilled charging, we will first charge the PEVs which depart at time stage $m+1$. Thus, if
\begin{equation}
\sum_{t=1}^{m} s_t -  \sum_{t=1}^{m} \sum_{ i \in \{ i|t_i^{(e)} = t \} } d_i \geq \sum_{ i \in \{ i|t_i^{(e)} = m+1 \} } d_i,
\end{equation}
we finish charging of PEVs with deadline $m+1$, and then go to charge the PEVs with deadline $m+2$.
If
\begin{equation}
\sum_{t=1}^{m} s_t -  \sum_{t=1}^{m} \sum_{ i \in \{ i|t_i^{(e)} = t \} } d_i < \sum_{ i \in \{ i|t_i^{(e)} = m+1 \} } d_i,
\end{equation}
then the PEVs with deadline $m+1$ have been charged as power $\sum_{t=1}^{m} s_t -  \sum_{t=1}^{m} \sum_{ i \in \{ i|t_i^{(e)} = t \} } d_i$.
At time stage $m+1$. Since
\begin{equation}\small
\sum_{t=1}^{m+1} s_t \geq \sum_{t=1}^{m+1} \sum_{ i \in \{ i|t_i^{(e)} = t \} } d_i,
\end{equation}
then,
\begin{equation} \small
s_{m+1} \geq \sum_{ i \in \{ i|t_i^{(e)} = m+1 \} } d_i - \left( \sum_{t=1}^{m} s_t -  \sum_{t=1}^{m} \sum_{ i \in \{ i|t_i^{(e)} = t \} } d_i\right),
\end{equation}
which means all the PEVs with deadline $m+1$ can be fulfilled charging. This is because we will charge the PEVs with deadline $m+1$ first by the EDF scheduling. Thus, there exists at least a set of $x_{it}$'s that satisfy
\begin{subequations}\label{theorem:eq:s3}
\begin{align}
&\sum_{t = t_i^{(s)}}^{t_i^{(e)}} x_{it} = d_i,  \forall i \in \{ i|t_i^{(e)} = m+1 \},\\
&x_{i,m+1} \geq 0,  \forall i \in \{ i|t_i^{(e)} = m+1 \}.
\end{align}
\end{subequations}
Combining (\ref{theorem:eq:s2}) and (\ref{theorem:eq:s3}), we get that all the PEVs whose deadline at or before stage $m+1$ can be fulfill charging, i.e., there exist  at least a set of $x_{it}$'s that satisfy
\begin{subequations}
\begin{align}
&\sum_{t = t_i^{(s)}}^{t_i^{(e)}} x_{it} = d_i,  \forall i \in \{ i|t_i^{(e)} \leq m+1 \},\\
&x_{it} \geq 0, \forall t = t_i^{(s)}, \cdots, t_i^{(e)},  \forall i \in \{ i|t_i^{(e)} \leq m+1 \}.
\end{align}
\end{subequations}
Therefore, we can conclude that by EDF scheduling, there always exists at least a set of $x_{it}$'s that is feasible to (\ref{model6}). This completes the proof.
\hfill $\blacksquare$

\subsection{Proof of Theorem \ref{theorem:convex}:}\label{appendix:proof:theorem:convex}
First, we show that if there exists a PEV parking in the station at both time $t_1$ and $t_2$, i.e.,
\begin{equation}
t_1, t_2 \in \{ t_i^{(s)}, \cdots, t_i^{(e)} \},
\end{equation}
and
\begin{equation}
x_{it_1}^* \geq 0, x_{it_2}^* > 0,
\end{equation}
then the optimal total loads at time $t_1$ and $t_2$ must satisfy that
\begin{equation}\label{eq:proof:theorem:convex}
s_{t_1}^* + l_{t_1} \geq s_{t_2}^* + l_{t_2}.
\end{equation}
The Karush-Kuhn-Tucker (KKT) conditions to the convex problem (\ref{model6}) are
\begin{subequations}\small
\begin{align}
\label{kkt1}
  &f'(\sum_{i \in \mathcal{I}(t)} x_{it} + l_t)
   - \lambda_i - \omega_{it} = 0,
   i \in \mathcal{N}, t = t_i^{(s)}, \cdots, t_i^{(e)}, \\
\label{kkt2}
  &\lambda_i(d_i - \sum_{t = t_i^{(s)}}^{t_i^{(e)}} x_{it}) = 0, i \in \mathcal{N}, \\
\label{kkt3}
  &\omega_{it}x_{it} = 0, i \in \mathcal{N}, t = t_i^{(s)}, \cdots, t_i^{(e)},
\end{align}
\end{subequations}
where $\lambda, \omega$ are the non-negative optimal Lagrangian multipliers corresponding to (\ref{model8}) and (\ref{model9}), respectively. We separate our analysis into the following two cases:
\begin{enumerate}
  \item If $x^*_{it_1} = 0$  for a particular PEV $i$ at a time slot $t_1 \in \{t_i^{(s)}, \cdots, t_i^{(e)}\}$, then, by complementary slackness, we have $ \omega_{it_1} > 0$. From (\ref{kkt1}),
      \begin{equation}\label{lemma:kkt:proof1}
      f'(s_{t_1} + l_{t_1}) = \lambda_i + \omega_{it_1}.
      \end{equation}
  \item If $x^*_{it_2} > 0$ for PEV $i$ during a time slot $t_2 \in \{t_i^{(s)}, \cdots, t_i^{(e)}\}$, we can infer from (\ref{kkt3}) that $\omega_{it_2} = 0$. Then,
      \begin{equation}\label{lemma:kkt:proof2}
      f'(s_{t_2}+ l_{t_2}) = \lambda_i.
      \end{equation}
\end{enumerate}
On the other hand, since $f(s_t + l_t)$ is a strictly convex function of $s_t + l_t$, then $f'(s_t + l_t)$ is an increasing function. From the above discussions, we get the following two conclusions:
\begin{enumerate}
  \item If $x^*_{it_1} > 0, x^*_{it_2} > 0$, then by (\ref{lemma:kkt:proof2}),
  \begin{equation}
  f'(s_{t_1} + l_{t_1}) = f'(s_{t_2} + l_{t_2} ) = \lambda_i.
  \end{equation}
  Due to the monotonicity of $f'(s_t)$, we have $s_{t_1}^* + l_{t_1} = s_{t_2}^* + l_{t_2}$.
  \item If $x^*_{it_1} = 0, x^*_{it_2} > 0$,  then by (\ref{lemma:kkt:proof1}) and (\ref{lemma:kkt:proof2}), there is
      \begin{equation}
      f'(s_{t_1} + l_{t_1} ) = \lambda_i + \omega_{it_1} > f'(s_{t_2} + l_{t_2}) = \lambda_i.
      \end{equation}
      Since $f'(s_t)$ is a increasing function, we have $s_{t_1}^* + l_{t_1} \geq s_{t_2}^* + l_{t_2}$.
\end{enumerate}

Consider two function $\hat{f}(s_t + l_t)$ and $\bar{f}(s_t + l_t)$. Let $\hat{x}_{it}^*$ and $ \bar{x}_{it}^*$ denote the optimal solutions to (\ref{model6}) with $f(s_t + l_t)$ replaced by $\hat{f}(s_t + l_t)$ and $\bar{f}(s_t + l_t)$, respectively.
Define $\hat{s}_t^*, \bar{s}_t^*$ as
\begin{equation}
\hat{s}_t^* = \sum_{i \in \mathcal{I}(t)} \hat{x}_{it}^*, \bar{s}_t^* = \sum_{i \in \mathcal{I}(t)} \bar{x}_{it}^*, t = 1, \cdots, T,
\end{equation}
respectively.
Suppose that there exists a time slot $t_1$ such that
\begin{equation}\label{theorem:convex:eq1}
\hat{s}_{t_1}^* < \bar{s}_{t_1}^*.
\end{equation}
Since
\begin{equation}
\sum_{t=1}^T \hat{s}_t^* = \sum_{t=1}^T \bar{s}_t^* = \sum_{i \in \mathcal{N}} d_i,
\end{equation}
there must exist another time slot $t_2$ such that
\begin{equation}\label{theorem:convex:eq2}
\hat{s}_{t_2}^* > \bar{s}_{t_2}^*
\end{equation}
and
\begin{equation}\label{theorem:convex:eq2}
\hat{s}_{t_1}^* + \hat{s}_{t_2}^* = \bar{s}_{t_1}^* + \bar{s}_{t_2}^*
\end{equation}
Thus, we can find a PEV $i \in \mathcal{N}$ such that
\begin{equation}
\hat{x}_{it_1}^* < \bar{x}_{t_1}^*, \hat{x}_{it_2}^* > \bar{x}_{it_2}^*.
\end{equation}
As a result,
\begin{equation}
\hat{x}_{it_2}^* > 0
\end{equation}
since $\bar{x}_{it_2}^* \geq 0$.
Based on (\ref{eq:proof:theorem:convex}), there is
\begin{equation}\label{theorem:convex:eq3}
\hat{s}_{t_2}^* + l_{t_2} \leq \hat{s}_{t_1}^* + l_{t_1}.
\end{equation}
Combining (\ref{theorem:convex:eq1})(\ref{theorem:convex:eq2})(\ref{theorem:convex:eq3}), we get
\begin{equation}\label{theorem:convex:eq4}
\bar{s}_{t_2}^* + l_{t_2} < \hat{s}_{t_2}^* + l_{t_2} \leq \hat{s}_{t_1}^* + l_{t_1}< \bar{s}_{t_1}^* + l_{t_1}.
\end{equation}
Since $\bar{f}(s_t+l_t)$ is a strictly convex function of $s_t+l_t$, then, based on \eqref{theorem:convex:eq2} and \eqref{theorem:convex:eq4}, we have
\begin{equation}\label{eq:theorem:convex:1}
\bar{f}(\bar{s}_{t_1}^* + l_{t_1}) + \bar{f}(\bar{s}_{t_2}^* + l_{t_2}) > \bar{f}(\hat{s}_{t_1}^* + l_{t_1}) + \bar{f}(\hat{s}_{t_2}^* + l_{t_2}).
\end{equation}
This contradicts with the fact that the $\bar{s}_t^*$ is the optimal total charging rate for objective function $\bar{f}(s_t + l_t )$. Therefore, the optimal charging solution $s_t^*$ is the same for any strictly convex function $f(s_t + l_t)$.
Next, we show that optimal solution $s_t^*$ is a load flattening solution that minimizes $\sum_{t =1}^T (s_t + l_t - \frac{ \sum_{t =1}^T s_t+l_t }{T})^2$ subjecting to \eqref{ws_modelb} and \eqref{ws_modelc}.
Based on the argument that $s_t^*$ is the same for any strictly convex function $f(s_t + l_t)$, then it is equivalent to show that $\sum_{t =1}^T (s_t + l_t - \frac{ \sum_{t =1}^T s_t+l_t }{T})^2$ is a strictly convex function of $s_t + l_t$.
Since
\begin{equation}
\frac{ \sum_{t =1}^T s_t+l_t }{T} = \frac{\sum_{i \in \mathcal{N}} d_i + \sum_{t =1}^T l_t}{T},
\end{equation}
which indicates that $\frac{ \sum_{t =1}^T s_t+l_t }{T}$ is a constant.
Then, we see that $\sum_{t =1}^T (s_t + l_t - \frac{ \sum_{t =1}^T s_t+l_t }{T})^2$ is a strictly convex function of $s_t + l_t$.
This completes the proof.
\hfill $\blacksquare$

\subsection{Proof of Proposition \ref{proposition:ev_lowerbound}: }
\label{appendix:proof:proposition:ev_lowerbound}
First, we show that $\Psi_{1}(\boldsymbol{\Xi})$ is a convex function of  $\boldsymbol{\Xi}$.
For any $\boldsymbol{\Xi}$, define $s_t^*(\boldsymbol{\Xi})$ as the optimal solution that minimizes $\Psi_{1}(\boldsymbol{\Xi})$ subject to (\ref{ws_modelb}) - (\ref{ws_modelc}).
For an arbitrary pair of $\Xi'$ and $\Xi''$,
let $\boldsymbol{\Xi}''' = \lambda \boldsymbol{\Xi}' + (1-\lambda) \boldsymbol{\Xi}'' \forall \lambda \in [0, 1].$
Then, the charging schedule $s_t(\boldsymbol{\Xi}''')$ such that $s_t(\boldsymbol{\Xi}''') = \lambda s_t^*(\boldsymbol{\Xi}') + (1-\lambda) s_t^*(\boldsymbol{\Xi}'')$
still satisfies (\ref{ws_modelb}) - (\ref{ws_modelc}) with $\boldsymbol{\Xi}=\boldsymbol{\Xi}'''$ due to the linearity of the constraints.
Based on the convexity of $f(s_t +l_t )$, we have
\begin{equation}\small\label{lemma:convex:proof1}
\begin{aligned}
&\sum_{t=1}^T f(s_t(\boldsymbol{\Xi}''') +l_t ) \\
\leq &\lambda \sum_{t=1}^T f(s_t^*(\boldsymbol{\Xi}') +l_t ) + (1 - \lambda) \sum_{t=1}^T f(s_t^*(\boldsymbol{\Xi}'') +l_t )
\end{aligned}
\end{equation}
for all $\lambda \in [0,1]$.
On the other hand, let $s_t^*(\boldsymbol{\Xi}''')$ be the optimal solution that minimizes $\sum_{t=1}^T f(s_t +l_t )$ subject to (\ref{ws_modelb}) - (\ref{ws_modelc}) with $\boldsymbol{\Xi} =\boldsymbol{\Xi}'''$. Then,
\begin{equation}\small\label{lemma:convex:proof3}
\sum_{t=1}^T f(s_t^*(\boldsymbol{\Xi}''') +l_t ) \leq \sum_{t =1}^T f(s_t(\boldsymbol{\Xi}''') +l_t ).
\end{equation}
Combining (\ref{lemma:convex:proof1}) and (\ref{lemma:convex:proof3}), we have
\begin{equation}\small
\begin{aligned}
&\sum_{t=1}^T f(s_t^*(\boldsymbol{\Xi}''') +l_t )  \\
\leq &\lambda \sum_{t =1}^T f(s_t^*(\boldsymbol{\Xi}') +l_t ) + (1 - \lambda) \sum_{t =1}^T f(s_t^*(\boldsymbol{\Xi}'') +l_t ).
\end{aligned}
\end{equation}
Thus, we have established the convexity of $\Psi_{1}(\boldsymbol{\Xi})$ over the set of $\boldsymbol{\Xi}$. 
Therefore,
we have
\begin{equation}\small
\mathrm E \left [\Psi_{1} ( \boldsymbol{\Xi}) \right ] \geq \Psi_{1} (\mathrm E[\boldsymbol{\Xi}]).
\end{equation}
On the other hand, based on the definition of $s_t, \tilde{d}_t^i, \mu_t^j, i = 1, \cdots, e_1, j = t, \cdots, e_t$, we have $\sum_{t =1}^{T} s_t = \sum_{t = 1}^{e_1}\tilde{d}_t^1 +  \sum_{t = 2}^T \sum_{j = t}^{e_t} \mu_t^j.$
Then, by Jensen's inequality \cite{boyd2004},
\begin{subequations}\small
\begin{align}
\Psi_{1} (\mathrm E[\boldsymbol{\Xi}])
&\geq  \sum_{t=1}^T f \left ( \frac{\sum_{t = 1}^{e_1}\tilde{d}_t^1 +  \sum_{t = 2}^T \sum_{j = t}^{e_t} \mu_t^j + \sum_{t=1}^T \nu_t}{T} \right ) \\
&= T f \left ( \frac{\sum_{t = 1}^{e_1}\tilde{d}_t^1 +  \sum_{t = 2}^T \sum_{j = t}^{e_t} \mu_t^j + \sum_{t=1}^T \nu_t}{T} \right ).
\end{align}
\end{subequations}
This completes the proof.
\hfill $\blacksquare$

\subsection{Proof of Proposition \ref{proposition:bound:onpoly}: }
\label{appendix:proof:proposition:bound:onpoly}
$s_t \leq \sum_{n=t}^T \tilde{d}_n^t$ holds for all stage $t$.
On the other hand, we have
\begin{equation}\small
\sum_{n=t}^T \tilde{d}_n^t 
\leq \sum_{n=t-1}^T \tilde{d}_n^{t-1} - \tilde{d}_{t-1}^{t-1} + \sum_{n =t}^{e_t} \eta_n^t
= \sum_{m \in \mathcal{M} } \sum_{n = t}^{e_m} \eta_n^m,
\end{equation}
where the inequality  holds since $s_{t-1} \geq \tilde{d}_{t-1}^{t-1}$ and $\mathcal{M}=\{m| e_m \geq t, m = 1, \cdots, t\}$
Thus, we have $s_t \leq \sum_{(m,n) \in \mathcal{O}(t) } \eta_n^m,$
where $\mathcal{O}(t)$ is a bounded set for $t = 1, \cdots, T$.
Therefore, $\mathrm E \left [ \sum_{t=1}^T f(\sum_{(m,n) \in \mathcal{O}(t) } \eta_n^m + \iota_t) \right ]$ is an upper bound of $\Phi_{3}$.
This completes the proof.
\hfill $\blacksquare$


\subsection{Proof of Lemma \ref{lemma:periodic}:}\label{appendix:proof:lemma:periodic}
We provide the proof by discussing the following two cases:

1)If $\bar{j} \geq \bar{i}+p$, which means that $[\bar{i}, \bar{j}]$ and $[\bar{i}+p, \bar{j}+p]$ overlaps with each other, then the density of interval $[\bar{i}, \bar{j}+p]$ is higher than that of $[\bar{i}, \bar{j}]$, and the density of $[\bar{i}, \bar{j}+2p]$ is higher than that of $[\bar{i}, \bar{j}+p]$. So on and so forth. Finally, we see that the interval $[\bar{i}, \bar{j}+(r-1)p]$ has the maximum density over region $\Pi_3$. Thus, we have $\hat{i}_3 = \bar{i}, \hat{j}_3 = \bar{j}+(r-1)p,$
      and
    \begin{equation}\small
    \hat{Z} = \frac{ \sum_{n=\bar{i}}^{\bar{j}+(r-1)p} (\sum_{m=n}^{k+e_{n}} \mu_{m}^{n} + \nu_{n}) }{ \bar{j}+(r-1)p-\bar{i} +1 }.
    \end{equation}
    Likewise, for the region $\{i, j| i = k, j = k+\hat{e}+1, \cdots, T \}$, we have $\hat{j}_2 = \bar{j}+(r-1)p,$
    and
    \begin{equation}\small
     \hat{Y} = \frac{ \sum_{n=k}^{\bar{j}+(r-1)p} (\sum_{m=n}^{k+e_{n}} \mu_{m}^{n} + \nu_{n}) }{ \bar{j}+(r-1)p-k +1 }.
    \end{equation}

2)If $\bar{j} < \bar{i}+p$,
then, the density of interval $[\bar{i}, \bar{j}]$ is higher than that of $[\bar{i}, \bar{j}+p]$, and the density of $[\bar{i}, \bar{j}+p]$ is higher than that of $[\bar{i}, \bar{j}+2p]$. So on and so forth. Finally, we see that the interval $[\bar{i}, \bar{j}]$ has the maximum density over region $\Pi_3$.
Thus, we have $\hat{i}_3 = \bar{i}, \hat{j}_3 = \bar{j},$
      and the corresponding maximum density
  \begin{equation}\small
    \hat{Z} = \frac{ \sum_{n=\bar{i}}^{\bar{j}} (\sum_{m=n}^{k+e_{n}} \mu_{m}^{n} + \nu_{n}) }{ \bar{j}-\bar{i} +1 }.
    \end{equation}
    For the region $\Pi_2$, if $\bar{j} \leq k+\hat{e}+1$, then 
    \begin{equation}\small
    \hat{Y} = \frac{ \sum_{n=k}^{k + \hat{e}+1} (\sum_{m=n}^{k+e_{n}} \mu_{m}^{n} + \nu_{n}) }{ \hat{e} +2 }.
    \end{equation}
    If $\bar{j} > k+\hat{e}+1$, then 
    \begin{equation}\small
    \hat{Y} = \frac{ \sum_{n=k}^{\bar{j}} (\sum_{m=n}^{k+e_{n}} \mu_{m}^{n} + \nu_{n}) }{ \bar{j} -k +1 }.
    \end{equation}

This completes the proof.
\hfill $\blacksquare$

\subsection{Proof of Lemma \ref{lemma:stationary}:}\label{appendix:proof:lemma:stationary}
Let $\rho(i, j)$ denote the maximum density of $[i, j]$.
For any $i = k, j=k+\bar{e}+1, \cdots, T$, the density of interval $[i, j]$ is given by
\begin{equation}
\begin{aligned}
\rho(i,j) &= \frac{ \sum_{k = 1}^{j-k} \sum_{t=1}^{j-2k+1} \mu_t + \sum_{t=k}^{k+\bar{e}} \tilde{d}_t^k + l_k + (j-k)\nu }{ j-k + 1 } \\
&= \frac{ \sum_{t=1}^{j-k}(j-k+1-t) \mu_t + \sum_{t=k}^{k+\bar{e}} \tilde{d}_t^k + l_k - \nu}{ j-k+1 } + \nu.
\end{aligned}
\end{equation}
To prove that the maximum density is achieved by setting $j =T$, we only need to show $\rho(i,j)$ is a non-decreasing function of $j$ for each given $i$, i.e.,
\begin{equation}\label{lemma:stationary:proof:3}
\rho(i,j) \leq \rho(i, j+1), \forall k+\bar{e}+1 \leq j \leq T-1.
\end{equation}
Since
\begin{equation}
\begin{aligned}
\frac{ \sum_{t=1}^{j-k}(j-k+1-t) \mu_t + \sum_{t=k}^{k+\bar{e}} \tilde{d}_t^k}{ j-k } \leq \sum_{t=1}^{j-k+1} \mu_t, \\
 k+\bar{e}+1 \leq j \leq T-1,
\end{aligned}
\end{equation}
we have
\begin{equation}\small
\begin{aligned}
&\rho(i,j+1) \\
= &\frac{ \sum_{t=1}^{j-k}(j-k+1-t) \mu_t + \sum_{t=k}^{k+\bar{e}} \tilde{d}_t^k + \sum_{t=1}^{j-k+1} \mu_t + l_k - \nu}{ j-k+1 } + \nu\\
\geq & \frac{ \sum_{t=1}^{j-k}(j-k+1-t) \mu_t + \sum_{t=k}^{k+\bar{e}} \tilde{d}_t^k + l_k - \nu}{ j-k } + \nu \\
= & \rho(i,j),
\end{aligned}
\end{equation}
which implies (\ref{lemma:stationary:proof:3}).
Hence, $Y$ is the maximum density of $[k, j], j = k+\bar{e}+1, \cdots, T$.
Next, we show that $Z$ is the maximum density of $[k+1, T]$.
For any $k+1 \leq i \leq j \leq T$, the density of interval $[i, j]$ is given by
\begin{equation}
\begin{aligned}
\rho(i,j) &= \frac{ \sum_{k = 1}^{j-i+1} \sum_{t=1}^{j-i+2-k} \mu_t  }{ j-i+1 } + \nu \\
&= \frac{ \sum_{t=1}^{j-i+1}(j-i+2-t) \mu_t }{ j-i+1 } + \nu.
\end{aligned}
\end{equation}
To prove that the maximum density is achieved by setting $i = k+1, j =T$, we only need to show $\rho(i,j)$ is a non-decreasing function of $j$ for each given $i$, i.e.,
\begin{equation}\label{lemma:stationary:proof:1}
\rho(i,j) \leq \rho(i, j+1), \forall k+1 \leq i \leq j \leq T-1,
\end{equation}
and a non-increasing function of $i$ for each given $j$, i.e.,
\begin{equation}\label{lemma:stationary:proof:2}
\rho(i,j) \geq \rho(i+1, j), \forall k+1 \leq i+1 \leq j \leq T.
\end{equation}
On one hand, since
\begin{equation}
\frac{ \sum_{t=1}^{j-i+1}(j-i+2-t) \mu_t }{ j-i+1 } \leq \sum_{t=1}^{j-i+2} \mu_t, \forall k+1 \leq i \leq j \leq T,
\end{equation}
we have
\begin{equation}\small
\begin{aligned}
\rho(i,j+1)
&= \frac{ \sum_{t=1}^{j-i+1}(j-i+2-t) \mu_t + \sum_{t=1}^{j-i+2} \mu_t}{ j-i+2 } + \nu \\
&\geq \frac{ \sum_{t=1}^{j-i+1}(j-i+2-t) \mu_t }{ j-i+1 }  + \nu\\
& = \rho(i,j),
\end{aligned}
\end{equation}
which implies (\ref{lemma:stationary:proof:1}).
On the other hand, as
\begin{equation}
\frac{ \sum_{t=1}^{j-i}(j-i+1-t) \mu_t }{ j-i } \leq \sum_{t=1}^{j-i+1} \mu_t, \forall k+1 \leq i \leq j \leq T,
\end{equation}
then
\begin{equation}
\begin{aligned}
\rho(i+1,j)
&= \frac{ \sum_{t=1}^{j-i}(j-i+1-t) \mu_t }{ j-i } + \nu\\
&\leq \frac{ \sum_{t=1}^{j-i}(j-i+1-t) \mu_t +  \sum_{t=1}^{j-i+1} \mu_t}{ j-i+1 } + \nu\\
& = \rho(i,j),
\end{aligned}
\end{equation}
which implies (\ref{lemma:stationary:proof:2}).
This completes the proof.
\hfill $\blacksquare$

\end{document}